\documentclass[12pt, oneside]{article}   	
\usepackage{geometry}                		
\usepackage{authblk}								
\usepackage{graphicx}
\usepackage{amssymb,amsmath,amsthm,bm}
\usepackage{mathtools,amsfonts}
\usepackage{epstopdf}
\usepackage{color}
\usepackage[numbers]{natbib}
\usepackage{hyperref}
\usepackage{setspace}
\usepackage{lineno}
\newcommand{\Var} {\mbox{$\rm{Var}$\,}}
\newcommand{\Prob} {\mbox{$\rm{Prob}$\,}}
\newtheorem{theorem}{Theorem}
\newtheorem{lemma}{Lemma}
\newtheorem{proposition}{Proposition}
\newtheorem{remark}{Remark}
\newtheorem{corollary}{Corollary}

\title{The stationary and quasi-stationary properties of neutral multi-type branching process diffusions}
\author[$\dagger$]{Conrad J.\ Burden}
\author[$\star$]{Robert C.\ Griffiths}

\affil[$\dagger$]{Mathematical Sciences Institute, Australian National University, Canberra, Australia 

		\texttt{conrad.burden@anu.edu.au}
		
		{\small ORCID id: 0000-0003-0015-319X, Researcher id: D-5556-2016 } } 
		
\affil[$\star$]{School of Mathematics, Monash University, Melbourne, Australia 

		\texttt{Bob.Griffiths@Monash.edu}
		
		{\small ORCID id: 0000-0001-7190-5104}
		} 

\date{}		

\begin{document}
\maketitle


\pagebreak
\begin{abstract}
The stationary asymptotic properties of the diffusion limit of a multi-type branching process with neutral mutations are studied.  
For the critical and subcritical processes the interesting limits are those of quasi-stationary distributions conditioned on non-extinction.  Pedagogical derivations 
are given for known results that the limiting distributions 
for supercritical and critical processes are found to collapse onto rays aligned with stationary eigenvectors of the mutation rate matrix, in agreement with 
discrete multi-type branching processes.  For the sub-critical process the previously unsolved quasi-stationary distribution is obtained to first order in 
the overall mutation rate, which is assumed to be small.  The sampling distribution over allele types for a sample of given finite size is found to agree 
to first order in mutation rates with the analogous sampling distribution for a Wright-Fisher diffusion with constant population size.  
\end{abstract}

\section{Introduction}
\label{intro}

A multi-type branching process, as defined in Chapter II of the book by \citet{Harris64}, describes discrete non-overlapping generations of a 
population which is partitioned into $d$ types.  In this paper we will assume $d$ to be finite.  
Individuals in the population at any time step give birth to a non-negative integer valued random 
number of offspring in the next generation.  The number of offspring per parent 
is identically and independently distributed across parents of a given type and across generations.  The expected number $m_{ij}$ of offspring 
of type-$j$ per parent of type-$i$ is assumed to be finite for all $i, j \in \{1, \ldots d\}$, and to allow for the possibility that a population 
can become extinct, at least one parental type has a non-zero probability of producing no offspring.  If $\rho$ is the maximal eigenvalue of 
the matrix $(m_{ij})_{i, j = 1}^d$, which we assume to be irreducible, the process is said to be subcritical, critical or supercritical according as $\rho < 1$, $\rho = 1$ or $\rho > 1$ 
respectively. 

As a population genetics model, a multi-type branching process has some similarities to haploid 
Wright-Fisher and Moran models, in which the total population size is usually assumed to be constant or to vary deterministically, 
rather than varying stochastically with time.  To see the similarities, decompose the expected number of offspring per parent as $m_{ij} = \lambda_i r_{ij}$, 
where the $r_{ij}$ are elements of a finite-state Markov transition matrix whose rows sum to 1.  The $r_{ij}$ are per-generation mutation rates between 
alleles, and the $\lambda_i$ carry information about the relative fitness of allele types~\citep{burden2018mutation}.  If the distribution of the number of 
offspring per parent is independent of parental type, and hence $\lambda_i = \lambda = \rho$, the model corresponds to neutral mutations.  In this case, 
the total population size, ignoring allele types, is effectively the $d = 1$ case and evolves as a Bienyam\'e-Galton-Watson (BGW) branching 
process~\citep{watson1875probability}.  

As is well known, the asymptotic probability of extinction at large times of a 
supercritical BGW process, $\rho > 1$, is equal to the stable fixed point of the moment generating function of number of offspring per parent, whereas 
for $\rho \le 1$ the population becomes asymptotically extinct with probability 1.  In those cases for which extinction is almost certain, 
the interesting asymptotic limit is the 
so-called quasi-stationary distribution of the population size conditioned on non-extinction.  For the critical case, $\rho = 1$, the weak asymptotic 
limit of the surviving population divided by the number of generations is exponentially distributed~\citep{yaglom1947certain}.  
For a review of quasi-stationary distributions for discrete-state models see van Doorn and Pollett~\cite{vanDoorn2013}, and for a review of 
continuous-state branching processes see Lambert~\cite{lambert2007quasi}.  

Known asymptotic 
results for discrete multi-type branching processes are listed in \citet[][p44]{Harris64} and \citet[][pp186-192]{Athreya_1972}.  In summary, suppose 
$\mathbf{Y}^*(\tau) = (Y_1^*(\tau), \ldots, Y_d^*(\tau))$ is the vector of population sizes of each type at time step $\tau$, conditioned on non-extinction.  
Provided certain conditions on the number of offspring per parent are met~\cite{kesten1966limit}, then the distributions of the scaled conditional population sizes 
$\mathbf{Y}^*(\tau)/\rho^\tau$ if $\rho > 1$, or $\mathbf{Y}^*(\tau)/\tau$ if $\rho = 1$, collapse onto a ray aligned with the stationary left eigenvector of 
$(m_{ij})_{i, j = 1}^d$ as $\tau \to \infty$.  
Moreover, for the critical multi-type branching process $\rho = 1$ the distribution along the ray is exponential.  If $\rho < 1$ and the second moments of the 
number of offspring per parent are finite~\cite{jivrina1962asymptotic}, the limiting distribution as $\tau \to \infty$ of $\mathbf{Y}^*(\tau)$ exists, 
is independent of the initial condition $\mathbf{Y}^*(\tau)$, has known first moments, and does not collapse onto a ray.   
\and \citet[][Theorem~1]{buiculescu1975quasi} has shown that the condition on the second moments of the number of offspring per parent can be considerably weakened.  

In this paper we are concerned with the asymptotic behaviour at large times of neutral multi-type branching processes in the diffusion limit.  
The diffusion limit of a 1-allele branching process, or equivalently 
a BGW process, was formulated and solved completely by \citet{feller1951diffusion}.  The diffusion limit of a multi-type branching process 
studied in this paper is a particular case of multi-type continuous state branching processes, which are characterised in 
\citet{li2010measure}, \citet{barczy2015stochastic} and \citet{caballero2017affine}.  
Our specific formulation is easily relatable to population genetics models, and can be found in \citet{burden2018mutation}.  Here, diffusion limit 
is meant in the sense of \citet{kimura1964diffusion}, where simultaneous limits are taken in which the continuum time 
between generations is taken to zero, the effective population size becomes infinite, and the per-generation mutation rate is taken to zero in such a 
way that mutation events along any lineage become a continuous-time Markov process with a finite rate matrix $\gamma$.  

In Section~\ref{sec:MultitypeBranching} the multi-type branching diffusion for $d$ types is introduced as the limit of a discrete multi-type branching process.  
Because this paper is restricted to neutral mutations, the marginal distribution of the total population size is equivalent to that of a  $d = 1$ branching process, 
also known as a Feller diffusion.   Section~\ref{sec:1type} is a summary of known results for Feller diffusions which will be needed for subsequent sections, 
paying particular attention to the asymptotic stationary limit. In particular, the quasi-stationary distribution of the surviving population has an exponential limit law 
in the subcritical case, and in the critical case is exponential provided the population size is scaled by the continuum time, consistent with the Yaglom limit.  

Rigorous results for the asymptotic stationary behaviours for critical and supercritical multi-type diffusions can be found in the continuous-state branching 
process literature~\citep{champagnat2008limit,Kyprianou18}.  In Section~\ref{sec:AsymptoticSuperAndCrit} we provide relatively straightforward derivations which should be accessible the mainstream 
population genetics community.  The quasi-stationary and stationary distributions respectively are seen to collapse onto rays aligned with the principle left 
eigenvectors of the mutation rate matrix, consistent with the known asymptotic limits of the discrete processes described above.  

The quasi-stationary limit of the subcritical multi-type diffusion is less straightforward, and is the subject of the main results of this paper.  
In Section~\ref{sec:QSSubcriticalMultiType} the quasi-stationary distribution is calculated to first order in 
an overall scaled mutation rate $\theta$ indicating the magnitude of the off-diagonal elements of diffusion limit mutation rate matrix.  
The precise definition of $\theta$ is given by Eq.~(\ref{smallRatesGamma}) below.    
The small-$\theta$ approximation is appropriate to many biologically realistic settings and has been applied to multi-allele Wright-Fisher 
diffusions~\citep{burden2016approximate,Burden_2018,Burden_2019}, to the mathematically equivalent boundary mutation model 
approximation to the Moran model~\citep{vogl2015inference,SchrempfHobolth17}, and to estimation of mutation rate parameters from site frequency 
data~\citep{vogl2014estimating,BurdenTang17,VoglMikulaBurden20}.  Higher order moments of the quasi-stationary distribution to first order in $\theta$ 
and sampling distributions are derived in Section~\ref{sec:HIgherOrderMoments}.  In Section~\ref{sec:ComparisonWithNumerical}, a numerical computation 
of the quasi-stationary distribution for the $d = 2$ sub-critical neutral branching process is compared with the approximate solution of 
Section~\ref{sec:QSSubcriticalMultiType} in order to gauge the range of validity of the small-$\theta$ approximation.  

Conclusions are drawn in Section~\ref{sec:Conclusions}.  

%
%

\section{Neutral multi-type branching diffusion for $d$ types}
\label{sec:MultitypeBranching}

Consider a BGW branching process with discrete generations $\tau = 0, 1, 2, \ldots$.  Assume the numbers of offspring per individual per 
generation are i.i.d. random variables, represented here by a generic random variable $S$ with $\Pr(S = 0) > 0$, 
$\mathbb{E}[S] =  \lambda, \Var(S) = \sigma^2$, 
with $\lambda$ and $\sigma^2$ finite.  If the total population size at time step $\tau$ is $Y(\tau)$, then $Y(\tau + 1) = \sum_{i = 1}^{Y(\tau)} S_i$.  
Suppose further that the population is divided into $d$ types with population counts $\mathbf{Y} = (Y_1, \ldots, Y_d)$, and that the probability 
of an offspring being of type-$j$ given their parent is of type-$i$ is $r_{ij}$, independently for each offspring.  Here $r_{ij} \ge 0$ and $\sum_{j = 1}^d r_{ij} = 1$.  
This is an example of a broader class of processes called multi-type branching process \citep[][Chapter~5]{Harris64,mode1971multitype,Athreya_1972}.  
More specifically, it corresponds to neutral mutations within 
a branching population, in the sense that the mean $\lambda$ and variance $\sigma^2$ of the number of offspring per parent are the same for all types.  

We note that a weaker requirement on the number of offspring per parent that $\mathbb{E}[S \log S] < \infty$ is necessary and sufficient for 
the limit theorems mentioned in the introduction to hold~\cite{buiculescu1975quasi,kesten1966limit}.   However, the stronger requirement of finite $\sigma^2$ will 
enable the diffusion limit as defined below, and is likely to be satisfied in practical applications to population genetics.  

The diffusion limit is obtained by defining a continuous time $t$ and scaled population $\widetilde{X}(t)$ by  
\begin{equation}	\label{diffLimit}
t = \frac{\sigma^2 \tau}{Y(0)}, \quad \widetilde{X}(t) = \frac{Y(\lfloor \tau\rfloor)}{Y(0)}, \quad \widetilde{\mathbf{X}}(t) = \frac{\mathbf{Y}(\lfloor \tau\rfloor)}{Y(0)},
\end{equation}
and by taking the limit 
$Y(0) \rightarrow \infty$, $\lambda \rightarrow 1$, $\sigma^2$ fixed and $r_{ij} \to 0$, in such a way that 
\begin{equation}	\label{alphaDef}
\alpha := \frac{Y(0) \log \lambda}{\sigma^2}, \quad \gamma_{ij} :=  \frac{Y(0)}{\sigma^2}  (r_{ij} - \delta_{ij}), \quad i, j = 1, \ldots, d, 
\end{equation}
remain fixed.  We take $\gamma$ to be an irreducible rate matrix.  
Note that $\alpha$ can be any real number, and that $(\gamma_{ij})_{i, j = 1}^d$ is an instantaneous rate matrix satisfying 
$\gamma_{ij} \ge 0$ for $i \ne j$ and $\sum_{j = 1}^d \gamma_{ij} = 0$.  
The resulting diffusion generator defined on $\mathbb{R}^d_+$ for $X(t)$, from the approximation $\widetilde {X}(t)$, is 
\begin{equation}		\label{scaledGeneratorMultitype}
{\cal L} = \frac{1}{2} \sum_{i = 1}^d x_i \frac{\partial^2}{\partial x_i^2} + \alpha \sum_{i = 1}^d x_i \frac{\partial}{\partial x_i}
				+ \sum_{i , j= 1}^d  \gamma_{ji} x_j \frac{\partial}{\partial x_i}.
\end{equation}
A detailed derivation of the forward Kolmogorov equation for the exponentially scaled population $\mathbf{Z}(t) = \mathbf{X}(t)/e^{\alpha t}$ 
is given in \citet[][Section~3]{burden2018mutation}.  Derivation of the forward equation for $\mathbf{X}(t)$ follows a similar path, and the result is 
easily seen to be consistent with this generator.  

Let $f(\bm{x})$ be a bounded continuous  function with second derivatives existing. Then a standard backward Kolmogorov equation is
\begin{equation}
\frac{d}{d t}\mathbb{E}_{X(0)}\big [f(\bm{X}(t))\big ] = \mathbb{E}_{X(0)}\big [{\cal L}f(\bm{X}(t))\big],
\label{query:01}
\end{equation}
where the right side is the expectation of the function $g$ defined by $g = {\cal L}f$.  
An elementary sketch of the derivation of Eq.~(\ref{query:01}) for a 1-dimensional diffusion process is in 
\cite{KT1981} p214. The multi-dimentional derivation follows in a similar style.  
 Define the Laplace transform, for $\phi_i > 0$, as 
\begin{equation}	\label{LTofXDefn}
\psi({\bm \phi}, t; \alpha, \mathbf{x}_0) = \mathbb{E} \left[ \left. e^{-\sum_{i = 1}^d\phi_i X_i(t)} \right| \mathbf{X}(0) = \mathbf{x}_0 \right],  
\end{equation}
where $\sum_{i = 1}^d x_{0i} = 1$.    
With ${\cal L}$ the generator (\ref{scaledGeneratorMultitype}) and $f(\bm{x}) = e^{-\sum_{i = 1}^d\phi_i X_i(t)}$, Eq.~(\ref{query:01}) leads to 
 \[
   \frac{\partial \psi}{\partial t}
 = \sum_{i=1}^d\left(-\frac{1}{2}\phi_i^2 + \alpha \phi_i + \sum_{j=1}^d \gamma_{ij}\phi_j \right) \frac{\partial \psi}{\partial \phi_i}.
 \]
The initial boundary condition is
$$
\psi({\bm \phi}, 0; \alpha, \mathbf{x}_0) = e^{-\sum_{i = 1}^d\phi_i x_{0i}(0)}.
$$

%
%

\section{Known results for $d = 1$ type}
\label{sec:1type}

For a neutral multi-type branching process, the generator of the total population $X = \sum_{i = 1}^d X_i$ is the $d = 1$ case of Eq.~(\ref{scaledGeneratorMultitype}).  
In this case the index $i$ and the final, $\gamma$-dependent, term in the generator no longer appear, and the initial condition is $X(0) = 1$.  
The solution~\citep{feller1951two} and its asymptotic properties~\citep{lambert2007quasi} are well known.  Here we summarise results which will be needed 
later in this paper.  

The Laplace transform $\psi_{\rm 1-allele}(\phi, t; \alpha) = \mathbb{E}[e^{-\phi X(t)}]$ is found by integrating along characteristic curves to be~\citep[][p236]{Cox78} 
\begin{eqnarray*} 
\psi_{\rm 1-allele}(\phi, t; \alpha) & = & \exp\left\{\frac{-\alpha\phi e^{\alpha t}}{\alpha + \frac{1}{2}(e^{\alpha t} - 1) \phi} \right\} \nonumber \\
			& = & \sum_{\ell = 0}^\infty e^{-\mu(t; \alpha)} \frac{\mu(t; \alpha)^\ell}{\ell!}(1 + \beta(t; \alpha)\phi)^{-\ell}, 
\end{eqnarray*}
where 
\begin{equation}	\label{muBetaDef}
\mu(t; \alpha) = \frac{2\alpha e^{\alpha t}}{e^{\alpha t} - 1}, \qquad \beta(t; \alpha) = \frac{e^{\alpha t} - 1}{2\alpha}.  
\end{equation}
We set $\mu(t; 0) = 2/t$ and $\beta(t; 0) =  t/2$.  

This is the Laplace transform of a point mass $e^{-\mu(t; \alpha)}$ at $x = 0$ representing the probability that the population becomes extinct at or before time $t$, 
plus a continuous Poisson-Gamma mixture for $x > 0$.  The resulting density is   
\begin{eqnarray}	\label{FellersPoissonGammaMixture}
f_X(x, t; \alpha) & = &  \delta(x) e^{-\mu(t; \alpha)}  \\
	& & +\, \sum_{\ell = 1}^\infty e^{-\mu(t; \alpha)} \frac{\mu(t; \alpha)^\ell}{\ell!} \frac{x^{\ell - 1}}{\beta(t; \alpha)^\ell (\ell - 1)!} e^{-x/\beta(t; \alpha)}, \quad x \ge 0,  \nonumber
\end{eqnarray}
where $\delta(x)$ is the Dirac delta function.  For the subcritical and critical cases, $\alpha \le 0$,  eventual extinction of the entire population 
is almost certain, and in the supercritical case, $\alpha > 0$, eventual extinction occurs with probability $e^{-2\alpha}$.  
Eq.~(\ref{FellersPoissonGammaMixture}) can be interpreted as a sum over the number of initial ancestral founders at $t = 0$, with 
$\mu(t; \alpha)$ the mean number of ancestral families surviving at time $t$, and each family size independently and exponentially distributed with mean $\beta(t; \alpha)$.  

Consider now the weak asymptotic limit of $X(t)$ as $t \to \infty$.  For the supercritical case the stationary limit is best understood in terms of the random variable 
$Z(t) = X(t) e^{-\alpha t}$ corresponding to the population size relative to the mean exponential growth.  
From Eq.~(\ref{muBetaDef}),   $\mu(t; \alpha) \to 2\alpha$ and $\beta(t; \alpha) e^{-\alpha t} \to 1/(2\alpha)$ as $t \to \infty$ for $\alpha > 0$, 
giving the asymptotic density of $Z(t)$ as  
\begin{eqnarray}	\label{supercritStat_fZ}
f_Z(z, \infty; \alpha) & = & \lim_{t \rightarrow\infty} e^{\alpha t} f_X(z e^{\alpha t} , t; \alpha) \nonumber \\
	& = & \delta(z) e^{-2\alpha} +  \sum_{\ell = 1}^\infty \frac{(2\alpha)^{2\ell}}{\ell! (\ell - 1)!} z^{\ell - 1} e^{-2\alpha (1 + z)}, \quad z \ge 0.  
\end{eqnarray}

For the subcritical and critical cases the interesting stationary 
limit as $t \rightarrow \infty$ is the quasi-stationary distribution corresponding to conditioning on survival of the population.  
The density corresponding to the random variable $X(t) | (X(t)>0)$ is 
\begin{equation}	\label{defG}
g_X(x, t; \alpha) = \frac{f_X(x, t; \alpha) - \delta(x) p_0(t)}{1 - p_0(t)}, \quad x \ge 0,    
\end{equation}
where $p_0(t) = e^{-\mu(t, \alpha)}$ is the survival probability.  The corresponding Laplace transform is 
\begin{eqnarray}	\label{zetaDef}
\zeta_{\rm 1-allele}(\phi, t; \alpha) & = & \frac{\psi_{\rm 1-allele}(\phi, t; \alpha) - p_0(t)}{1 - p_0(t)} \nonumber \\
		& = &   \sum_{\ell = 1}^\infty \frac{e^{-\mu(t; \alpha)}}{1 - e^{-\mu(t; \alpha)}} \frac{\mu(t; \alpha)^\ell}{\ell!}(1 + \beta(t; \alpha)\phi)^{-\ell}.
\end{eqnarray}
For $\alpha < 0$ we have $\mu \to 0$ and $\beta \to 1/(2|\alpha|)$ as $t \to \infty$.  Only the $\ell = 1$ term in Eq.~(\ref{zetaDef}) survives the limit, leading to 
\begin{equation}	\label{1AlleleSubcritZeta}
\zeta_{\rm 1-allele}(\phi, \infty; \alpha) = \frac{2|\alpha|}{\phi + 2|\alpha|}, 	\qquad \alpha < 0,   
\end{equation}
which is the Laplace transform of the quasi-stationary exponential distribution
\begin{equation} \label{gExponential}
g_X(x, \infty; \alpha) = 2|\alpha| e^{-2|\alpha| x}, 	\qquad x \ge 0, \alpha < 0.   
\end{equation}

For the critical case, define the random variable   
\begin{equation*}
W(t) = \frac{X(t)}{t}.    
\end{equation*}
The density function conditioned on non-extinction of $W(t)|(W(t) > 0)$, is 
\begin{equation*} 
g_W(w, t) = tg_X(tw, t; 0),   
\end{equation*}
where the function $g$ is defined in Eq.~(\ref{defG}).  The corresponding Laplace transform is 
\begin{eqnarray*}
\zeta_W(\phi, t) & = & \zeta(\phi t^{-1}, t; 0) \nonumber \\
	& = & \sum_{\ell = 1}^\infty \frac{e^{-\mu(t; 0)}}{1 - e^{-\mu(t; 0)}} \frac{\mu(t; 0)^\ell}{\ell!}\left(1 + \beta(t; 0)\frac{\phi}{t}\right)^{-\ell}.
\end{eqnarray*}	
Once again only the $\ell = 1$ term survives the limit, giving 
\begin{equation}	\label{zetaStationary}
\zeta_W(\phi, \infty) = \frac{1}{1 + \tfrac{1}{2}\phi}, \qquad g_W(w, \infty) = 2e^{-2w}, \qquad \alpha = 0. 
\end{equation}
Reinstating the original variables for the discrete BGW process via Eq.~(\ref{diffLimit}) gives 
$
\lim_{\tau \to \infty} \Prob(Y(\tau)/\tau > z \mid Y(\tau) > 0) = e^{-2z/\sigma^2}. 
$
This agrees with Yaglom's well-known exponential limit law of \cite{yaglom1947certain}, a proof of which appears in \citet[][p20]{Athreya_1972}.   
Note that asymptotically, the entire surviving population is descended from a single ancestor from the initial population at time $\tau = 0$ 
in both the critical and sub-critical cases.    

%
%

\section{Asymptotic behaviours of neutral multi-type branching diffusions: supercritical and critical cases}
\label{sec:AsymptoticSuperAndCrit}

The main purpose of this paper is to study the asymptotic stationary behaviour of neutral multi-type diffusion processes conditional on non-extinction. 
The supercritical and critical cases 
are essentially covered in the existing literature using other methods and will be dealt with first.  The subcritical case is less straightforward and will be covered in subsequent 
sections.  

The supercritical case has previously been studied in detail by \citet{burden2018mutation}, with emphasis on the $d = 2$ case, and a more formal 
treatment in terms of measure valued processes for any finite number of types is to be found in \citet[][Theorem~1.4]{Kyprianou18}.  
Here we provide a derivation of the asymptotic stationary distribution for $d$ types by adapting and generalising the $d = 2$ proof in \citet[][Section~6.1]{burden2018mutation}.  

%
%
\begin{proposition}\label{proposition1}
Define the exponentially scaled variable 
\begin{equation*}
\mathbf{Z}(t) = \left . \mathbf{X}(t) e^{-\alpha t} \right| (\mathbf{X}(0) = \mathbf{x}_0).  
\end{equation*}
If $\alpha > 0$, the limit stationary distribution of $\mathbf{Z}(t)$ as $t\to \infty$ has a density 
\begin{equation}	\label{supercritMultiSoln}
f_\mathbf{Z}(\mathbf{z}, \infty; \alpha, \mathbf{x}_0) = \pi_d^{d - 2} f_Z(z, \infty; \alpha) \prod_{\ell = 1}^{d - 1} \delta(\pi_d z_\ell - \pi_\ell z_d), 
\end{equation}
where $z = \sum_{i = 1}^d z_i$, $\bm \pi = (\pi_1 \cdots \pi_d)$ is the stationary left eigenvector of $(\gamma_{ij})_{i, j = 1}^d$, and the function 
$f_Z(\cdot)$ is the single allele density for the total population defined by Eq.~(\ref{supercritStat_fZ}).  
\end{proposition}
\begin{proof}
The generator of $Z(t)$ acting on bounded continuous functions $g(\bm{z})$ with second derivatives existing is
\begin{equation}	\label{genOfZ}
{\cal L}_t = \frac{1}{2} e^{-\alpha t}\sum_{i = 1}^d z_i \frac{\partial^2}{\partial z_i^2} 
				+ \sum_{i , j= 1}^d  \gamma_{ji} z_j \frac{\partial}{\partial z_i}.
\end{equation}
The limit of ${\cal L}_t$ as $t\to \infty$ is
\begin{equation}
{\cal L}_\infty = \sum_{i , j= 1}^d  \gamma_{ji} z_j \frac{\partial}{\partial z_i}.
\label{newL:0}
\end{equation}
A stationary limit distribution is defined as one where $\mathbb{E}\big [{\cal L}_\infty g(\bm{Z})\big ]=0$ for functions in the domain of ${\cal L}_\infty$. Choosing $g(\bm{z}) = e^{-\sum_{i=1}^d\phi_iz_i}$ and denoting the Laplace transform 
\begin{equation*}	
\psi_\mathbf{Z}({\bm \phi}, t; \alpha, \mathbf{x}_0) = \mathbb{E} \left[e^{-\sum_{i = 1}^d\phi_i Z_i(t)} \right],  
\end{equation*} 
the stationary equation for the Laplace transform is
\begin{equation}	\label{asympPsiZpde}
\sum_{i,j=1}^d \gamma_{ji}\phi_i\frac{\partial \psi_\mathbf{Z}({\bm \phi}, \infty; \alpha, \mathbf{x}_0)}{\partial \phi_j} = 0. 
\end{equation}
A boundary condition is determined by setting ${\bm \phi} = (\phi, \cdots, \phi)$ and noting that the total population size $X = \sum_{i = 1}^d X_i$  evolves as 
a Feller diffusion for 1 allele type.  Thus  
\begin{equation}	\label{asympPsiZbc}
\psi_\mathbf{Z}(\phi{\mathbf 1}, \infty; \alpha, \mathbf{x}_0)) = \mathbb{E}\left[e^{-\phi Z}\right] = \psi_Z(\phi, \infty), 
\end{equation}
where $\psi_Z(\phi, \infty)$ is the Laplace transform of Eq.~(\ref{supercritStat_fZ}).  

In general, the irreducible rate matrix $\gamma$ has a complete set of left eigenvectors $\mathbf{v}^{(\ell)}$ and right eigenvectors $\mathbf{u}^{(\ell)}$ with 
normalisation condition $\mathbf{v}^{(k)} \cdot \mathbf{u}^{(\ell)} = \delta_{k\ell}$, and corresponding 
eigenvalues $\nu_\ell$, where $\ell = 0, \ldots d - 1$.  Specifically, $\mathbf{v}^{(0)} = {\bm \pi}$ is the left stationary eigenvector, 
$\mathbf{u}^{(0)} = \mathbf 1$, and $\nu_0 = 0$.  Suppressing the $\alpha$ dependence to simplify the notation, without loss of generality set 
\begin{equation*}
\psi_\mathbf{Z}({\bm \phi}, \infty) = h(\bm \phi \cdot \mathbf{v}^{(0)}, \ldots, \bm \phi \cdot \mathbf{v}^{(d - 1)}), 
\end{equation*}
where the function $h$ is to be determined.  Then Eq.(\ref{asympPsiZpde}) becomes 
\begin{equation*}
\sum_{\ell = 1}^{d - 1} \nu_\ell \bm\phi \cdot \mathbf{v}^{(\ell)} \partial_\ell h(\bm \phi \cdot \mathbf{v}^{(0)}, \ldots, \bm \phi \cdot \mathbf{v}^{(d - 1)}) = 0, 
\end{equation*}
where $\partial_\ell$ means partial differentiation with respect to the $\ell$th argument, $\ell = 0, ..., d - 1$.  Because this differential equation does not involve 
$\bm \phi \cdot \mathbf{v}^{(0)} = \bm \phi \cdot \bm \pi$, the function $h$ factors into 
\begin{equation}	\label{psiZFactored}
\psi_\mathbf{Z}({\bm \phi}, \infty) = h_0(\bm \phi \cdot \bm \pi) h_\perp(\bm \phi \cdot \mathbf{v}^{(1)} \ldots, \bm \phi \cdot \mathbf{v}^{(d - 1)}), 
\end{equation}
where $h_\perp$ satisfies 
\begin{equation*}
\sum_{\ell = 1}^{d - 1} \nu_\ell \xi_\ell \frac{\partial h_\perp(\xi_1, \ldots, \xi_{d - 1})}{\partial \xi_\ell} = 0. 
\end{equation*}
The characteristic curves parametrised by $s$, say, for this first order equation are determined from $d\xi_\ell/ds = \nu_\ell \xi_\ell$, for  $\ell = 1, \ldots, d - 1$. 
By solving these ordinary differential equations and eliminating $s$, it is easy to see that the characteristics can be stated as 
\begin{equation*}
\xi_1^{\nu_1} = c_\ell \xi_\ell^{\nu_\ell}, \qquad \ell = 2, \ldots, d - 1, 
\end{equation*}
with the set of constants $\{c_2, \ldots, c_{d - 1}\}$ labelling a characteristic.  Each characteristic passes through the origin and there is a characteristic curve passing through every 
point in the space spanned by $(\xi_1, \ldots, \xi_{d - 1})$.  Thus $h_\perp(\xi_1, \ldots, \xi_{d - 1}) = h_\perp(0, \ldots, 0)$ is constant throughout 
its domain provided $h_\perp(0, \ldots, 0)$ is well defined and finite.  

From the boundary condition Eq.~(\ref{asympPsiZbc}), we have that 
\begin{equation*}
h_0(\phi)h_\perp(0, \ldots, 0) = \psi_Z(\phi, \infty), 
\end{equation*}
since for $\ell > 0$,  $\phi{\mathbf 1} \cdot \mathbf{v}^{(\ell)} = \phi\mathbf{u}^{(0)} \cdot \mathbf{v}^{(\ell)} = 0$ by the orthogonality condition.  Returning to 
Eq.~(\ref{psiZFactored}), we have 
\begin{equation*}
\psi_\mathbf{Z}({\bm \phi}, \infty) = h_0(\bm \phi \cdot \bm \pi) h_\perp(0, \ldots, 0) = \psi_Z(\bm \phi \cdot \bm \pi, \infty),  
\end{equation*}
where $\psi_Z(\phi, \infty)$ is the Laplace transform of the 1-allele solution $f_Z(z, \infty; \alpha)$. 
It is straightforward to check that this is the Laplace transform of Eq.~(\ref{supercritMultiSoln}).  
\end{proof}

The interpretation is that the distribution collapses onto a line density aligned with the stationary eigenvector of the rate matrix, and, conditional on 
the population not becoming extinct, the proportion $Z_i/Z$ of allele type-$i$ in the population converges almost surely to $\pi_i$.  This result is the continuum 
version of a particular case of the limit theorem for a discrete supercritical BGW process stated in \citet[][Theorem~9.2]{Harris64} or 
\citet[][p19]{mode1971multitype}.  
Numerical computations for $d = 2$ and very small mutation rates by Burden and Wei~\citep{burden2018mutation} have displayed a collapse of the distribution onto a line density 
that begins after a rapid changeover point at $\alpha t \approx - \log\theta$, where $\theta$ is a measure of the overall mutation rate (see Eq.(\ref{smallRatesGamma}) below).  
These computations showed that the dynamics was dominated by the exponentially scaled genetic drift term in Eq.~(\ref{genOfZ}) before the changeover point, 
and and by the mutation term after the changeover point \citep[see Eq.~(57) of ][]{burden2018mutation}.  Heuristically, one expects a similar rapid changeover in the limit of small mutation rates for general $d$.  

The critical case also exhibits a collapse onto a line density aligned with the stationary eigenvector of the rate matrix in the asymptotic limit, except that, 
because extinction of the population is almost certain, the appropriate limit is the quasi-stationary distribution.  
%
%
\begin{proposition}\label{proposition2}
Define the scaled random variable
\begin{equation*}
\mathbf{W}(t) = \frac{\mathbf{X}(t)}{t}.
\end{equation*}
Then if $\alpha = 0$, the limit stationary distribution of the scaled population conditioned on non-extinction, 
$\mathbf{W}(t) \mid (|\mathbf{W}(t)| >0)$ as $t \to \infty$, is 
\begin{equation}	\label{critMultiSoln}
g_\mathbf{W}(\mathbf{w}, \infty) = 2e^{-2w} \pi_d^{d - 2} \prod_{\ell = 1}^{d - 1} \delta(\pi_d w_\ell - \pi_\ell w_d). 
\end{equation}
\end{proposition}
\begin{proof}
We first determine the (unconditional) generator for $\mathbf{W}(t)$ by considering a Laplace transform argument.
Let $\mathcal{L}$ be as in Eq.~(\ref{scaledGeneratorMultitype}) 
with $\alpha = 0$, then
\begin{equation*}
\begin{aligned}
&\frac{d}{d t} \mathbb{E}\Big [e^{-\sum_{i=1}^d\phi_i W_i(t)}\Big ]\\
&~~= \frac{\partial}{\partial t} \mathbb{E}\Big [e^{-\sum_{i=1}^d\phi_i X_i(t)/t}\Big ]\\
&~~= t^{-2}\mathbb{E}\Big [\sum_{i=1}^d\phi_iX_i(t)e^{-\sum_{i=1}^d\phi_i X_i(t)/t}\Big ] + 
\mathbb{E}\Big [{\cal L}e^{-\sum_{i=1}^d\phi_iX_i(t)/t}\Big ]\\
&~~= t^{-2}\mathbb{E}\Big [\sum_{i=1}^d\phi_iX_i(t)e^{-\sum_{i=1}^d\phi_i X_i(t)/t}\Big ]\\
 &~~~~+\mathbb{E}\Big [\Big (\frac{1}{2}t^{-2}\sum_{i=1}^d\phi_i^2X_i(t) - t^{-1}\sum_{i,j=1}^d\gamma_{ji}X_j(t)\phi_i\Big )
 e^{-\sum_{i=1}^d\phi_iX_i(t)/t}\Big ]\\
&~~=  -\frac{1}{t}\mathbb{E}\Big [\sum_{i=1}^dW_i(t)\frac{\partial}{\partial W_i(t)}e^{-\sum_{i=1}^d\phi_i W_i(t)}\Big ]\\
 &~~~~+\mathbb{E}\Big [\Big (\frac{1}{2t}\sum_{i=1}^dW_i(t)\frac{\partial^2}{\partial W^2_i(t)} + \sum_{i,j=1}^d\gamma_{ji}W_j(t)\frac{\partial}{\partial W_i(t)}\Big )
 e^{-\sum_{i=1}^d\phi_iW_i(t)}\Big ].\\
\end{aligned}
\end{equation*}
The generator for $\mathbf{W}(t)$ is therefore 
\begin{equation}	\label{unconditionedWGen}
\mathcal{L}_\mathbf{W} = \frac{1}{2t} \sum_{i = 1}^d w_i \frac{\partial^2}{\partial w_i^2} 
					+ \sum_{i, j = 1}^d \gamma_{ji} w_j \frac{\partial}{\partial w_i} 
					- \frac{1}{t} \sum_{i = 1}^d w_i \frac{\partial}{\partial w_i}.  
\end{equation}

Let $\mathbb{E}^*$ be expectation in the distribution conditional on non-extinction.  Let the generator in the conditional distribution be ${\cal L}^*_t$. Suppose $g$ is a bounded continuous function with second derivatives. Then
 \[
 \mathbb{E}^*[g(\bm{W}(t))] = \frac{\mathbb{E}\big [g(\bm{W}(t))]  - g( \bm{0} )p_0(t)}{1-p_0(t)}
 \]
 and
 \begin{eqnarray*}
 \frac{d}{d t}\mathbb{E}^*[g(\bm{W}(t))]
 	& = & \frac{\mathbb{E}\big [{\cal L}_{\bm{W}}g]}{1-p_0(t)} 
 		+ \frac{p_0^\prime(t)}{1-p_0(t)}\frac{\mathbb{E}[g(\bm{W}(t))- g(\bm{0})]}{1-p_0(t)}\\
 	& = & \frac{\mathbb{E}\big [{\cal L}_{\bm{W}}g]}{1-p_0(t)} 
		 + \frac{p_0^\prime(t)}{1-p_0(t)}\frac{\mathbb{E}[g(\bm{W}(t))- g(\bm{0})p_0(t)]}{1-p_0(t)} 
		 - g(\bm{0}) \frac{p_0^\prime(t)}{1-p_0(t)}\\
  	& =  & \frac{\mathbb{E}\big [{\cal L}_{\bm{W}}g]}{1-p_0(t)} 
 		+ \frac{p_0^\prime(t)}{1-p_0(t)}\mathbb{E}^*[g(\bm{W}(t))]
 		- g(\bm{0}) \frac{p_0^\prime(t)}{1-p_0(t)}.
 \end{eqnarray*}
 That is
 \begin{equation}
 {\cal L}^*_tg = {\cal L}_{\bm{W}}g + \frac{p_0^\prime(t)}{1-p_0(t)}(g - g(\bm{0})).
 \label{condgen:250}
 \end{equation}
 With $p_0(t) = e^{-\mu(t;0)} = e^{-2/t}$,  
 \[
 \frac{p_0^\prime(t)}{1-p_0(t)} \approx 1/t \to 0 \text{~as~} t \to \infty,
  \]
 so the limit generator is 
 \begin{equation}
 {\cal L}^*_\infty =  \sum_{i,j=1}^d \gamma_{ji}w_i\frac{\partial}{\partial w_j}.
 \label{limitgen:0}
 \end{equation}

A stationary limit distribution is defined as one where $\mathbb{E}^*\big [{\cal L}^*_\infty g(\bm{W})\big ]=0$ for functions in the domain of ${\cal L}^*_\infty$. 
Choosing $g(\bm{w}) = e^{-\sum_{i=1}^d\phi_iw_i}$ and denoting the Laplace transform 
\begin{equation*}	
\zeta_\mathbf{W}({\bm \phi}, t) = \mathbb{E}^* \left[e^{-\sum_{i = 1}^d\phi_i W_i(t)} \right],  
\end{equation*} 
the stationary equation for the Laplace transform is
\begin{equation*}
 \sum_{i, j = 1}^d \gamma_{ji}\phi_i\frac{\partial \zeta_\mathbf{W}({\bm \phi}, \infty)}{\partial \phi_j}  = 0.   
 \end{equation*}
A boundary condition is determined by setting ${\bm \phi} = (\phi, \cdots, \phi)$ and noting that for neutral mutations the total scaled population size 
$W = \sum_{i = 1}^d W_i$  evolves as the case for 1 allele type.  Thus by Eq.~(\ref{zetaStationary}),  
\begin{equation*}
\zeta_\mathbf{W}(\phi{\mathbf 1}, \infty) = \lim_{t \to \infty}\mathbb{E}\left[\left. e^{-\phi W(t)} \right| W(t) > 0 \right] = (1 + \tfrac{1}{2}\phi)^{-1}.  
\end{equation*}
The method of solution is identical to that for the asymptotic supercritical case, and leads to 
\begin{equation*}
\zeta_\mathbf{W}({\bm \phi}, \infty) = \left(1 + \tfrac{1}{2} {\bm \phi}\cdot{\bm \pi} \right)^{-1}, 
\end{equation*}
the inverse Laplace transform of which is Eq.~(\ref{critMultiSoln}).  
\end{proof}

Note that in the 1-dimensional case the Laplace transform of $X(t)/t$ tends to zero, so even though (\ref{limitgen:0}) is also the limit from the unconditioned 
generator (\ref{unconditionedWGen}) it does not give the correct solution because $X(t)/t$ does not have a finite limit and the solution is that the Laplace transform is zero. 
In the supercritical 1-dimensional case $e^{-\alpha t}X(t)$ converges to a proper limit, so it is not necessary to condition on survival.

The interpretation of Proposition~\ref{proposition2} is that the distribution collapses onto a line density of magnitude $2e^{-2w}$ aligned with the stationary eigenvector of the rate matrix $\gamma$.  
In other words, conditional on the population not becoming extinct, the proportion $X_i/X= W_i/W$ of allele type-$i$ in the population converges almost surely to $\pi_i$.  
This result is the diffusion limit analogue of \citet[][Theorem~1, p191]{Athreya_1972}.

%
%

\section{Quasi-stationary limit of a subcritical multi-type branching diffusion}
\label{sec:QSSubcriticalMultiType}

A complete solution of the quasi-stationary density for the subcritical case remains intractable.  In the following we derive an approximation to the 
quasi-stationary density which is correct to first order in small mutation rates.  We begin with two lemmas.  

%
%
\begin{lemma}\label{lemma1}
Define the Laplace transform of the multi-type population $\mathbf{X}(t)$ conditioned on survival of the population as 
$$
\zeta({\bm \phi}, t; \alpha, \mathbf{x}_0) = \mathbb{E} \left[ \left. e^{-\sum_{i = 1}^d\phi_i X_i(t)} \right| \mathbf{X}(0) = \mathbf{x}_0, X(t) > 0 \right], 
$$
where $X(t) = \sum_{i = 1}^d X_i(t)$.  Then if $\alpha < 0$, the Laplace transform of the limiting quasi-stationary distribution, 
 $\zeta({\bm \phi}) \equiv \zeta(\bm{\phi}, \infty; \alpha, \mathbf{x}_0)$, satisfies 
   \begin{eqnarray}
 \sum_{i=1}^d\left(-\frac{1}{2}\phi_i^2 - |\alpha| \phi_i + \sum_{j=1}^d \gamma_{ij}\phi_j\right)
 \frac{\partial \zeta}{\partial \phi_i} - |\alpha|(1-\zeta)=0.
 \label{LT:50}
 \end{eqnarray}
\end{lemma}
\begin{proof}
Let $\mathbb{E}^*$ be expectation in the distribution conditional on non-extinction, so that for any bounded continuous function $g$ with 
second derivatives,
\[
 \mathbb{E}^*[g(\bm{X}(t))] = \frac{\mathbb{E}\big [g(\bm{X}(t))]  - g( \bm{0} )p_0(t)}{1-p_0(t)}.
\]
Following the same argument as that leading to Eq.~(\ref{condgen:250}), the generator in the conditional distribution acting on $g$ is 
\[
 {\cal L}^*_tg = {\cal L}g + \frac{p_0^\prime(t)}{1-p_0(t)}(g - g(\bm{0})),
 \]
where ${\cal L}$ is given by Eq.~(\ref{scaledGeneratorMultitype}).  From Eq.~(\ref{muBetaDef}),
 \[
 \lim_{t\to \infty}\frac{p_0^\prime(t)}{1-p_0(t)} = -\alpha, \qquad \alpha < 0,  
 \]
 and thus the limit generator acting on $g$ is 
\[
 {\cal L}^*_\infty g = \left(\frac{1}{2} \sum_{i = 1}^d x_i \frac{\partial^2}{\partial x_i^2} - |\alpha| \sum_{i = 1}^d x_i \frac{\partial}{\partial x_i}
				+ \sum_{i , j= 1}^d  \gamma_{ji} x_j \frac{\partial}{\partial x_i}\right) g - |\alpha| (g - g(\bm{0})).
\]
Choosing $g(\mathbf{x}) = e^{-\sum_{i = 1}^d\phi_i x_i(t)}$ and setting $\mathbb{E}^*[{\cal L}^*_\infty g(\mathbf{X})] = 0$ then leads to Eq.~(\ref{LT:50}). 
 \end{proof}
%
%
\begin{lemma}\label{lemma2}
For a subcritical multi-type process, the mean of the limiting quasi-stationary distribution conditional on survival of the population is 
\begin{equation}	\label{mu_i}
\mu_i  := \left. - \frac{\partial \zeta({\bm \phi})}{\partial \phi_i} \right|_{{\bm \phi} = {\bm 0}} = \frac{1}{2|\alpha|} \pi_i, \qquad \alpha < 0,\ i = 1, \ldots, d, 
\end{equation}
where $\zeta({\bm \phi})$ is the solution to Eq.~(\ref{LT:50}), and 
$\bm{\pi}$ is the left stationary eigenvector of the rate matrix $\gamma$, normalised so that $\sum_{i = 1}^d \pi_i = 1$.  
\end{lemma}
\begin{proof}
Differentiating Eq.~(\ref{LT:50}) with respect to $\phi_r$ gives 
\begin{equation*}
-\phi_r \frac{\partial \zeta}{\partial \phi_r} + \sum_{i = 1}^d \gamma_{ir} \frac{\partial \zeta}{\partial \phi_i}
			+ \sum_{i = 1}^d \left(-\frac{1}{2}\phi_i^2 - |\alpha| \phi_i + \sum_{j = 1}^d \gamma_{ij} \phi_j \right) \frac{\partial^2 \zeta}{\partial \phi_i \partial \phi_r} = 0, 
\end{equation*}
and setting $\bm{\phi} = \bm{0}$ then gives 
\begin{equation*}
\sum_{i = 1}^d \gamma_{ir} \mu_i = 0. 
\end{equation*}
Furthermore, setting $\bm{\phi} = \bm{1} \phi = (1, \ldots, 1) \phi$ for scalar $\phi$ in Eq.~(\ref{LT:50}), and noting that Eq.~(\ref{LTofXDefn}) and 
the first line of Eq.~(\ref{zetaDef}) 
imply $\zeta(\bm{1} \phi) = \zeta_{\rm{1-allele}}(\phi)$ and that the chain rules implies $\partial \zeta(\bm{1} \phi)/\partial\phi_i = - \mu_i + \mathcal{O}(\phi)$ 
as $\phi \to 0$, gives 
\begin{equation*}
(\tfrac{1}{2}\phi^2 + |\alpha| \phi) \left(\sum_{i  = 1}^d\mu_i + \mathcal{O}(\phi)\right)  - |\alpha|(1-\zeta_{\rm{1-allele}}(\phi))=0.  
\end{equation*}
Substituting from Eq.~(\ref{1AlleleSubcritZeta}), dividing through by $\phi$ and then setting $\phi = 0$ then gives 
\begin{equation} \label{normalisationMu}
\sum_{i = 1}^d \mu_i = \frac{1}{2|\alpha|}. 
\end{equation}
Thus $\bm{\mu}$ is the stationary left eigenvalue of $\gamma$, normalised by Eq.~(\ref{normalisationMu}), as required.  
\end{proof}

%
%

\subsection{Small mutation rates}
\label{sec:SmallRates}
When studying small mutation rates, a convenient parameterisation for the rate matrix is 
\begin{equation}	\label{smallRatesGamma}
\gamma_{ij} = \tfrac{1}{2} \theta (P_{ij} - \delta_{ij}), 
\end{equation}
where $P_{ij} \ge 0$ are the elements of a finite state Markov transition matrix satisfying $\sum_{j = 1}^d P_{ij} =1$.  As noted in \citet{Burden_2019}, 
$\theta$ is arbitrary up to the constraint
$$
\tfrac{1}{2} \theta \ge \max_{i = 1, \ldots, d} \sum_{j: j \ne i} \gamma_{ij} = \max_{i = 1, \ldots, d} (-\gamma_{ii}), 
$$
and the choice of $\theta$ determines the $P_{ij}$.  Specifically, for a parent-independent rate matrix (PIM) satisfying $\gamma_{ij} = \gamma_j$ (independent of $i$) 
for $i \ne j$, the canonical parameterisation is $\tfrac{1}{2}\theta = \sum_{i = 1}^d \gamma_i$, which ensures $P_{ij} = \pi_j$, where $(\pi_1 \cdots \pi_d)$ is the 
stationary left-eigenvector of the rate matrix.  

Now consider a subcritical multi-type branching diffusion with general small mutation rates as in Eq.~(\ref{smallRatesGamma}) where $\theta << 1$.  
The differential equation Eq.~(\ref{LT:50}) for the Laplace transform of the subcritical quasi-stationary distribution density is scale invariant, and without loss of 
generality one can set $\alpha = -\frac{1}{2}$ to obtain
 \begin{equation} \label{LTAlphaEquals1}
 \sum_{i=1}^d\left[-\phi_i(1 + \phi_i) 
		+   \theta \sum_{j=1}^d (P_{ij} - \delta_{ij}) \phi_j \right] \frac{\partial \zeta}{\partial \phi_i} - (1-\zeta) = 0.
 \end{equation}
Results for any $\alpha < 0$ can be reconstructed by making replacements $\phi \to \frac{1}{2}|\alpha|^{-1}\phi$, 
$\theta \to \frac{1}{2}|\alpha|^{-1}\theta$ and $\zeta \to \zeta$, and the corresponding quasi-stationary density can be reconstructed from 
\begin{equation}	\label{gScaling}
g_\mathbf{X}(\mathbf{x}, \infty; \alpha, \theta) = 2|\alpha| g_\mathbf{X}\left(2|\alpha|\mathbf{x}, \infty; -\tfrac{1}{2}, \tfrac{1}{2}|\alpha|^{-1} \theta\right).  
\end{equation}

%
%
\begin{theorem}\label{theorem1}
The first order in $\theta$ solution to Eq.~(\ref{LTAlphaEquals1}) is 
\begin{equation}	\label{zetaExpansion}
\zeta(\bm \phi) = \zeta_0(\bm \phi) + \theta \zeta_1(\bm \phi) + o(\theta), 
\end{equation}
as $\theta \to 0$, where 
\begin{equation}	\label{zeta0}
\zeta_0(\bm \phi) = \sum_{i = 1}^d \pi_i \left(1 + \phi_i \right)^{-1}, 
\end{equation} 
and 
\begin{equation}	\label{zeta1}
\zeta_1(\bm \phi) = -\sum_{i, j = 1}^d \pi_j P_{ji} 
		(1 + \phi_i) \left( \left(1 + \phi_j\right)^{-1} - \left(1 + \phi_i \right)^{-1} \right)^2 \phi_i^{-1} \log \left(1 + \phi_i \right).  
\end{equation}
\end{theorem}
\begin{proof}
When $\theta = 0$ the $d$ types decouple, so $\zeta_0$ must be a linear combination of 1-allele solutions of the form of Eq.~(\ref{1AlleleSubcritZeta}).  
Furthermore, for agreement with the $\theta \to 0$ limit, the first moments must be as in Eq.~(\ref{mu_i}), and thus $\zeta_0$ is as given in Eq.~(\ref{zeta0}).  

Now work on the second term $\zeta_1$.  Assuming Eq.~(\ref{zetaExpansion}) and equating the coefficient of $\theta$ in Eq.~(\ref{LTAlphaEquals1}), 
\begin{eqnarray}	\label{zeta1PDE}
- \sum_{i = 1}^d \phi_i \left(1 + \phi_i  \right) \frac{\partial \zeta_1 }{\partial \phi_i} & = & 
	-\zeta_1 - \sum_{i, j=1}^d (P_{ij} - \delta_{ij}) \phi_j \frac{\partial \zeta_0}{\partial \phi_i} \nonumber \\
	& = & -\zeta_1 + \sum_{i, j=1}^d (P_{ij} - \delta_{ij}) \phi_j \frac{\pi_i}{\left(1 + \phi_i \right)^{2}}.
\end{eqnarray}
This equation is solved by integrating $\zeta_1$ along characteristic curves parametrised by a parameter $s$, say, in $\bm\phi$ space.  
These curves satisfy  
\begin{equation*}
\frac{d\phi_i}{ds} = -\phi_i \left(1 + \phi_i  \right).
\end{equation*}
For each $i = 1, \ldots, d$, 
\begin{equation}	\label{sMinusCi}
s - c_i = \log \left| \frac{1 + \phi_i}{\phi_i} \right|, 
\end{equation}
with $c_1, \ldots, c_d$ integration constants.  
It suffices to restrict $\bm\phi$ to the positive sector, giving 
\begin{equation}	\label{phiChars}
\phi_i = \frac{e^{-(s - c_i)}}{1 - e^{-(s - c_i)}}, \qquad \phi_i > 0, 
\end{equation}
as plotted in Fig.~\ref{fig:CharacteristicPlot}(a).  
For the characteristic passing through a given point $\bm{\phi}$, the integration constants $c_i$ are determined up to an overall additive constant 
independent of $i$ by 
\begin{equation*}
 c_j - c_i = \log \left| \frac{(1 + \phi_i)\phi_j}{(1 + \phi_j)\phi_i} \right|.
\end{equation*}
Arbitrarily choosing any one of the $c_i$ determines the remaining $d - 1$ integration constants.  The one-parameter family of 
characteristics for $d = 2$ are plotted in Fig.~\ref{fig:CharacteristicPlot}(b). 

\begin{figure}[t!]
\begin{center}
\centerline{\includegraphics[width=0.9\textwidth]{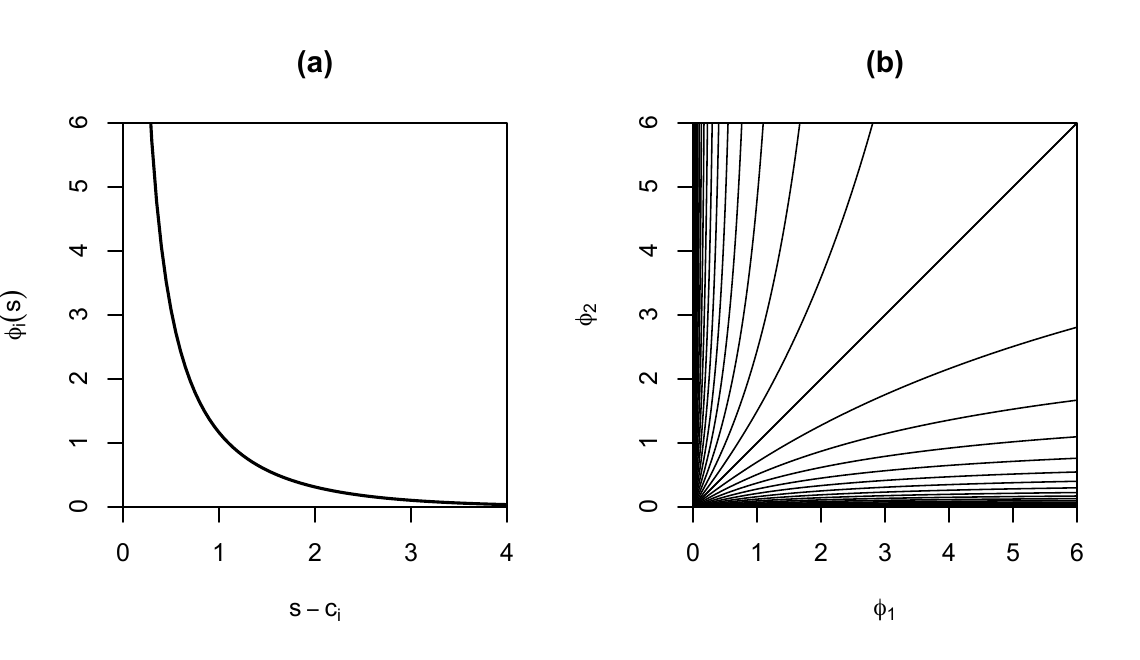}}
\caption{(a) Solutions $\phi_i(s)$ to the characteristic equations for the Laplace transform of the quasi-stationary, subcritical multitype branching diffusion.  
(b) Characteristic curves for $d = 2$ types.  } 
\label{fig:CharacteristicPlot}
\end{center}
\end{figure}

Along the characteristic passing through any given $\bm\phi$,  Eq.~(\ref{zeta1PDE}) implies  
\begin{equation*}
\frac{d\zeta_1}{ds} + \zeta_1  =  \sum_{i, j=1}^d \pi_i P_{ij} \frac{\phi_j}{(1 + \phi_i)^2} - \sum_{i = 1}^d \frac{\pi_i \phi_i}{(1 + \phi_i)^2}. 
\end{equation*}
Substituting Eq.~(\ref{phiChars}) and multiplying through by the integrating factor $e^s$, gives 
\begin{equation*} 
\frac{d}{ds}\left(\zeta_1 e^s\right) =  \sum_{i, j=1}^d \pi_j P_{ji} e^{c_i} \frac{\left(1 - e^{c_j - s}\right)^2}{1 - e^{c_i - s}} 
								- \sum_{i=1}^d \pi_i e^{c_i} \left(1 - e^{c_i - s}\right).  
\end{equation*}
For the integral of the first term, we need 
\begin{eqnarray*}
I_{ij}(s, \bm{c}) & = & \int \frac{\left(1 - e^{c_j - s}\right)^2}{1 - e^{c_i - s}} ds \nonumber \\
		& = & e^{2(c_j - c_i)}\left(s + e^{c_i - s}\right) + 2s e^{c_j - c_i} \left(1 - e^{c_j - c_i}\right)  \nonumber \\
		& & \qquad\qquad + \left(1 - e^{c_j - c_i}\right)^2 \left( s + \log \left|1 - e^{c_i - s} \right|\right), 
\end{eqnarray*}
up to an arbitrary constant which may depend on $\bm c$.  From Fig.~\ref{fig:CharacteristicPlot}(a) it is clear that $\max_i(c_i) < s < \infty$ for $\phi_i > 0$, 
so the absolute value signs in the last line can be dispensed with.  Then 
\begin{equation*}
\zeta_1 = \sum_{i, j=1}^d \pi_j P_{ji} e^{c_i - s} I_{ij}(s, \bm{c}) - \sum_{i=1}^d \pi_i e^{c_i - s} \left(s + e^{c_i - s}\right) + K(\bm{c}) e^{-s}, 
\end{equation*}
where $K(\bm{c})$ is a characteristic-dependent integration constant.  
It is straightforward to check by making use of the fact that $\sum_{j = 1}^d \pi_j P_{ji} = \pi_i$ that the terms proportional to $s$ cancel, leaving 
\begin{eqnarray*}
\zeta_1 & = &  \sum_{i, j=1}^d \pi_j P_{ji} e^{c_i - s} \left\{ e^{2(c_j - c_i) + c_i - s} + \left(1 - e^{c_j - c_i}\right)^2 \log \left(1 - e^{c_i - s} \right) \right\} \nonumber \\
		& & \qquad\qquad -  \sum_{i=1}^d \pi_i e^{2(c_i - s)}   + K(\bm{c}) e^{-s} \nonumber \\
		& = & \sum_{i, j=1}^d \pi_j P_{ji} e^{c_i - s} \left(1 - e^{c_j - c_i}\right)^2 \log \left(1 - e^{c_i - s} \right) + K(\bm{c}) e^{-s}, 
\end{eqnarray*}
where terms have been cancelled in the last line by making use of $\sum_{i = 1}^d P_{ji} = 1$.  

Since $\zeta_1$ is determined by a specified point $\bm{\phi}$ through which the characteristic passes, Eq.~(\ref{sMinusCi}) implies that $\zeta_1$ must depend on 
$s$ and $c_1, \ldots, c_d$ only via combinations of $s - c_i$.  Thus 
\begin{equation*}
K(\bm{c}) e^{-s} = \sum_{i = 1}^d b_i e^{c_i - s}, 
\end{equation*} 
for constants $b_1, \ldots, b_d$.  Reinstating the $\phi_i$ via Eq.~(\ref{sMinusCi}) then gives 
\begin{eqnarray*} 
\zeta_1(\bm{\phi}) & = & \sum_{i, j = 1}^d \pi_j P_{ji} \frac{\phi_i}{1 + \phi_i} \left(1 - \frac{\phi_j(1 + \phi_i)}{\phi_i(1 + \phi_j)}\right)^2 \log \left(\frac{1}{1 + \phi_i} \right) 
						+ \sum_{i = 1}^d b_i \frac{\phi_i}{1 + \phi_i} 	\nonumber \\
			& = & \sum_{i, j =1}^d \pi_j P_{ji}  \frac{1 + \phi_i}{\phi_i} \left( \frac{1}{1 + \phi_j} - \frac{1}{1 + \phi_i} \right)^2 \log \left(\frac{1}{1 + \phi_i} \right) 
								 + \sum_{i = 1}^d b_i \frac{\phi_i}{1 + \phi_i}. 
\end{eqnarray*}
The $b_i$ are determined from the first moments.  Expanding in powers of $\phi_i$, 
\begin{equation*}
\zeta_1(\bm{\phi}) = - \sum_{i, j =1}^d \pi_j P_{ji} (\phi_i - \phi_j)^2 (1 + {\cal O}(\phi)) + \sum_{i = 1}^d b_i \phi_i (1 + {\cal O}(\phi)).  
\end{equation*}
Only the second term contributes to the first moment, leading to 
\begin{equation*}
\mu_i = - \left. \frac{\partial}{\partial \phi_i} (\zeta_0 + \theta \zeta_1 + o(\theta)) \right|_{\bm{\phi} = 0} = \pi_i - \theta b_i + o(\theta). 
\end{equation*}
Comparing with the exact result to all orders in $\theta$, Eq.~(\ref{mu_i}) with $|\alpha| = \tfrac{1}{2}$, we see that the first moments are accounted for 
by $\zeta_0$, and thus $b_i = 0$, giving $\zeta_1$ as in Eq.~(\ref{zeta1}).  
\end{proof}
\begin{remark}
The inverse Laplace transform of $\zeta_0$ is 
\begin{equation}	\label{g0}
g_0(\mathbf{x}) = \sum_{i = 1}^d \pi_i e^{-x_i} \prod_{j \ne i} \delta(x_j),  
\end{equation}
which represents an exponentially distributed  line density along each $x_i$-axis.  
\end{remark}

The following lemma is needed before inverting the Laplace transform to $\mathcal{O}(\theta)$.  
%
%
\begin{lemma}\label{lemma3}
For any real $a > 0$, the Laplace transform of 
\begin{equation} \label{invLTIn3rdTerm}
- x^{a\theta - 1} E_2(x) + \left(\frac{1}{a\theta} - \gamma \right) \delta(x), 
\end{equation}
is 
\begin{equation}	\label{aDependentLT}
\left(1 + \phi\right) \phi^{-1} \log \left(1 + \phi \right) + \mathcal{O}(\theta), 
\end{equation}
as $\theta \to 0$, where  \citep[][Eq.~5.1.4]{Abramowitz:1965sf}
\begin{equation*}
E_n(z) = \int_1^\infty \frac{e^{-zt}}{t^n} dt, 
\end{equation*}
is the exponential integral.
\end{lemma}
\begin{proof}
Expanding the logarithm, 
\begin{eqnarray}	\label{expandedLT3rdTerm}
\left(1 + \phi\right) \phi^{-1} \log \left(1 + \phi \right) 
	& = & - \left(1 + \phi\right) \sum_{n = 1}^\infty \frac{(-1)^n}{n} \phi^{n - 1} \nonumber \\
	& = & - \sum_{n = 0}^\infty \frac{(-1)^{n + 1}}{n + 1} \phi^n -  \sum_{n = 1}^\infty \frac{(-1)^n}{n} \phi^n \nonumber \\
	& = & 1 - \sum_{n = 1}^\infty \frac{(-1)^n}{n(n + 1)}  \phi^n.  
\end{eqnarray}
We are required to check that this agrees with the Laplace transform of Eq.~(\ref{invLTIn3rdTerm}) to $\mathcal{O}(\theta)$.  First note that 
for any $\alpha > 0$, 
\begin{eqnarray*}
\int_0^\infty x^{\alpha - 1} E_2(x) \, dx & = & \int_0^\infty x^{\alpha - 1} \int_1^\infty \frac{e^{-xt}}{t^2} dt dx \nonumber \\
	& = & \int_1^\infty \frac{1}{t^2} \left( \int_0^\infty x^{\alpha - 1} e^{-xt} dx \right) dt \nonumber \\
	& = &\Gamma(\alpha) \int_1^\infty \frac{dt}{t^{2 + \alpha}} \nonumber \\
	& = & \frac{\Gamma(\alpha)}{1 + \alpha}.  
\end{eqnarray*}
Then the Laplace transform of the first term in Eq.~(\ref{invLTIn3rdTerm}) is 
\begin{eqnarray} \label{invLT1stTermInThirdTerm}
- \lefteqn{\int_0^\infty x^{a\theta - 1} E_2(x) e^{-x\phi} dx } & & \nonumber \\
	& = & - \sum_{n=0}^\infty \frac{(-1)^n}{n!} \phi^n \int_0^\infty x^{a\theta + n - 1} E_2(x) dx \nonumber \\
	& = & - \frac{\Gamma(a\theta)}{1 + a\theta} - 
				\sum_{n=1}^\infty \frac{(-1)^n}{n!} \frac{\Gamma(n)}{1 + n}\phi^n \left(1 + \mathcal{O}(\theta)\right)  \nonumber \\
	& = & -\frac{1}{a\theta} + 1 + \gamma  - \sum_{n=1}^\infty \frac{(-1)^n}{n(n + 1)} \phi^n + \mathcal{O}(\theta),
\end{eqnarray}
where we have used $\Gamma(z) = z^{-1} - \gamma + \mathcal{O}(z)$ as $z \to 0$ in the last line.  The Laplace transform of the second term 
in Eq.~(\ref{invLTIn3rdTerm}) is $1/(a\theta) - \gamma$, which, when added to Eq.~(\ref{invLT1stTermInThirdTerm}) agrees with 
Eq.~(\ref{expandedLT3rdTerm}) up to $\mathcal{O}(\theta)$.  
\end{proof}
In the following proofs we use the notation $(i \leftrightarrow j)$ to mean an expression where $i$ and $j$ are exchanged in an immediately preceding expression.
%
%
\begin{theorem}\label{theorem2}
The inverse Laplace transform of Eq.~(\ref{zetaExpansion}) is 
\begin{equation}	\label{gToOrderTheta}
g_{\mathbf{X}}(\mathbf{x}) =  \sum_{1 \le i < j \le d} g^{\rm surface} (x_i, x_j) \prod_{\ell \ne i, j} \delta(x_\ell) + 
											 \sum_{i = 1}^d g^{\rm line} (x_i) \prod_{\ell \ne i} \delta(x_\ell)  + o(\theta), 
\end{equation}
where 
\begin{equation} 	\label{gSurface}
g^{\rm surface} (x_i, x_j) = \theta \pi_j P_{ji} \left\{ x_j e^{-x_j} x_i^{a_{ij}(x_j)\theta - 1} E_2(x_i) + 
										2e^{-x_j} E_1(x_i) \right\} + (i \leftrightarrow j), 
\end{equation}
is a density over the 2-dimensional surface spanned by the $x_i$ and $x_j$ axes,  
\begin{equation} 	\label{gLine} 
g^{\rm line} (x_i) = - \theta \pi_i(1 - P_{ii})  \left\{E_1(x_i) + [\gamma(1 - x_i) + \log x_i] e^{-x_i} \right\},  
\end{equation}
is a line density along the $x_i$-axis, and  $a_{ij}(x_j)$, $i, j = 1, \ldots, d$, $i \ne j$, are a set of functions constrained by 
\begin{equation}	\label{choiceOfA}
\sum_{i \ne j} P_{ji} x_j \frac{1}{a_{ij}(x_j)} = 1, \qquad j = 1, \ldots, d.
\end{equation}
Eq.~(\ref{gToOrderTheta}) is the required first order in $\theta$ density of the quasi-stationary distribution for a subcritical multi-type 
branching diffusion.  
\end{theorem}
\begin{proof}
Expanding Eq.~(\ref{zeta1}), 
\begin{eqnarray}\label{expandedZeta1}
\zeta_1(\bm{\phi}) & = & -\sum_{(i, j):i \ne j} \pi_j P_{ji} \left\{ \left(1 + \phi_i\right)^{-1}  \phi_i^{-1} \log \left(1 + \phi_i \right) \right.\nonumber \\
	& & \qquad\qquad                 -2            \left(1 + \phi_j\right)^{-1}  \phi_i^{-1} \log \left(1 + \phi_i \right) \nonumber \\ \nonumber \\
	& & \qquad\qquad                  +\,  \left. \left(1 + \phi_j\right)^{-2}  \left(1 + \phi_i\right) \phi_i^{-1} \log \left(1 + \phi_i \right) \right\} \nonumber \\
	& = & -\sum_{(i, j):i \ne j} \pi_j P_{ji} \{\widetilde T_1 + \widetilde T_2 + \widetilde  T_3\},  
\end{eqnarray} 
say.  

The first two terms are inverted by making use of the results that the Laplace transform of $e^{-x}$ is $\left(1 + \phi\right)^{-1}$ and the 
Laplace transform of the exponential integral $E_1(x)$ is $ \phi^{-1} \log \left(1 + \phi \right)$.  
The inverse transform of $\widetilde T_1$ is the convolution integral  
\begin{eqnarray*}
\int_0^{x_i} e^{-(x_i - u)} E_1(u)\, du & = & e^{-x_i } \int_0^{x_i} e^u \int_1^\infty \frac{e^{-ut}}{t} \,dt du \nonumber \\
	& = &  e^{-x_i } \int_1^\infty \frac{1}{t} \frac{1 - e^{-x_i(t - 1)}}{t - 1} \, dt  \nonumber \\
	& = & e^{-x_i } \mathcal{I}(x_i), 
\end{eqnarray*}
where 
\begin{equation*}
\mathcal{I}(x) = \int_1^\infty \frac{1 - e^{-x(t - 1)}}{t(t - 1)} \, dt.
\end{equation*}
We have that $\mathcal{I}(0) = 0$ and 
\begin{equation*}
\mathcal{I}'(x) = e^x \int_1^\infty \frac{e^{-xt}}{t} dt = e^xE_1(x),  
\end{equation*}
and thus 
\begin{eqnarray*}
\mathcal{I}(x) & = & \lim_{\epsilon \to 0} \int_\epsilon^x e^u E_1(u) du \nonumber \\
	& = & \lim_{\epsilon \to 0} \left[ \left. e^u E_1(u) \right|_{u=\epsilon}^x + \int_\epsilon^x \frac{du}{u} \right] \nonumber \\
	& = & e^x E_1(x) + \gamma + \log x, 
\end{eqnarray*}
where $\gamma$ is the Euler-Mascheroni constant, and we have used that \citep[][Eq.~5.1.11]{Abramowitz:1965sf} 
$E_1(z) = -\gamma - \log z + \mathcal{O}(z)$ as $z \to 0$.  Thus the inverse Laplace transform of $\widetilde T_1$ is 
\begin{equation}	\label{invLT1stTerm}
T_1 = \left\{E_1(x_i) + [\gamma + \log x_i] e^{-x_i }\right\} \prod_{\ell \ne i} \delta(x_\ell).  
\end{equation}

For $i \ne j$, the inverse Laplace transform of $\widetilde T_2$ is 
\begin{equation}	\label{invLT2ndTerm}
T_2 = -2e^{-x_j} E_1(x_i) \prod_{\ell \ne i,j} \delta(x_\ell). 
\end{equation}

Inverting $\widetilde T_3$ requires Lemma~\ref{lemma3} for the $x_i$-dependent factors, and that the Laplace transform of 
$\left(1 + \phi_j\right)^{-2}$ is $x_j e^{-x_j}$ for the $x_j$-dependent factor.    
Furthermore, by carrying out the Laplace transform first as an integral over $x_i$, and then as 
an integral over $x_j$, it is clear that any dependence of the introduced parameter $a$ on $x_j$ can be absorbed into the $\mathcal{O}(\theta)$ part of 
Eq.~(\ref{aDependentLT}).  Thus the inverse Laplace transform of $\widetilde T_3$ is 
\begin{equation}	\label{invLT3rdTerm}
T_3 = - x_j e^{-x_j} \left[x_i^{a_{ij}(x_j)\theta - 1} E_2(x_i) - \left(\frac{1}{a_{ij}(x_j)\theta} - \gamma \right) \delta(x_i) \right] \prod_{\ell \ne i, j} \delta(x_\ell). 
\end{equation}

Reassembling the parts from Eqs. (\ref{g0}), (\ref{expandedZeta1}), (\ref{invLT1stTerm}), (\ref{invLT2ndTerm}) and (\ref{invLT3rdTerm}), 
the inverse Laplace transform of Eq.(\ref{zetaExpansion}) is 
\begin{eqnarray}	\label{gReassembled}
g_{\mathbf{X}}(\mathbf{x}) & = & \sum_{j = 1}^d \pi_j e^{-x_j} \prod_{\ell \ne j} \delta(x_\ell) 
				 - \theta \sum_{(i, j): i\ne j} \pi_j P_{ji} (T_1 + T_2 +  T_3) + o(\theta) \nonumber \\
	& = & \sum_{j = 1}^d \pi_j  e^{-x_j} \left( 1 - 
								\sum_{i \ne j} P_{ji} x_j \frac{1}{a_{ij}(x_j)} \right) \prod_{\ell \ne j} \delta(x_\ell) \nonumber \\
	&  & -\, \theta \sum_{(i, j): i\ne j} \pi_j P_{ji}   \left\{\left(E_1(x_i) + [\gamma + \log x_i] e^{-x_i }\right) \delta(x_j) 
			- 2e^{-x_j} E_1(x_i) \vphantom{x^()}\right. \nonumber \\ \nonumber \\
	&  &  \qquad \left. - \, x_j e^{-x_j} \left[x_i^{a_{ij}(x_j)\theta - 1} E_2(x_i) + \gamma \delta(x_i) \right] \right\} \prod_{\ell \ne i, j} \delta(x_\ell)
				 + o(\theta) .    \nonumber \\
\end{eqnarray}
Recall that the zero-th order solution, Eq.~(\ref{g0}), is a set of line densities representing the $\theta \to 0$ limit of singular behaviour near each axis.  By choosing 
$a_{ij}(x_j)$ to satisfy Eq.~(\ref{choiceOfA}),  
the leading order term is removed and the singular behaviour near each axis is exposed in a term containing a factor $x_i^{a_{ij}(x_j)\theta - 1}$ 
arising from the final line of Eq.~(\ref{gReassembled}).  The resulting density becomes  
\begin{eqnarray*}
g_{\mathbf{X}}(\mathbf{x}) & = & \theta \sum_{(i, j): i\ne j} \pi_j P_{ji} \left\{ x_j e^{-x_j} x_i^{a_{ij}(x_j)\theta - 1} E_2(x_i) + 
										2e^{-x_j} E_1(x_i) \right\} \prod_{\ell \ne i, j} \delta(x_\ell) \nonumber \\
	& &  \qquad -\, \theta \sum_{i = 1}^d \sum_{j \ne i} \left\{ \pi_j P_{ji}  \left(E_1(x_i) + [\gamma + \log x_i] e^{-x_i }\right) \right.   \nonumber \\
	& &         \qquad\qquad\qquad \left. -\pi_i P_{ij} \gamma x_i e^{-x_i}\right\} \prod_{\ell \ne i} \delta(x_\ell)  + o(\theta) ,
\end{eqnarray*}
which is equivalent to Eqs.(\ref{gToOrderTheta}), (\ref{gSurface}) and (\ref{gLine}).  
\end{proof}
\begin{remark}
Note that $g^{\rm surface}$ and $g^{\rm line}$ are both invariant with respect to the arbitrary choice of $\theta$ in Eq.~(\ref{smallRatesGamma}).  
\end{remark}
\begin{remark}
We have not explicitly calculated the functions $a_{ij}(x_j)$ occurring in the surface density, except to state the constraint Eq.~(\ref{choiceOfA}).  
These functions serve the purpose of ensuring that singular behaviour of $g^{\rm surface}$ near the boundary of the positive $(x_i, x_j)$ quadrant 
remains integrable and that the density 
is correctly normalised.  For the purpose of calculating higher order moments of $X_i/X$ to $\mathcal{O}(\theta)$, and hence sampling distributions, it will 
turn out that in general the functions $a_{ij}(x_j)$  can be set to zero, that is, the behaviour $x_i^{a\theta - 1}$ can simply be replaced by $x_i^{-1}$.  
\end{remark}
%
%

\section{Higher order moments of the subcritical quasi-stationary distribution}
\label{sec:HIgherOrderMoments}

%
\subsection{Moments in $\mathbf{X}$  to order $\theta$}

In the following theorems $H_n = \sum_{k = 1}^n k^{-1}$ is the $n$th harmonic number for $n \ge 1$, and $H_0 := 0$.  
%
%
\begin{theorem}\label{theorem3}
Define moments in $\mathbf{X}$ for the quasi-stationary distribution by 
\begin{equation*}
\mathbb{E}_{\rm qs}\left[ \prod_{i = 1}^d X_i^{n_i} \right]  = 
	\int_0^\infty \cdots \int_0^\infty \left( \prod_{i = 1}^d x_i^{n_i} \right) g_{\mathbf{X}}(\mathbf{x}) d^d\mathbf{x}.  
\end{equation*}
Then for integer $n \ge 1$, and $r \in \{1, \ldots, d\}$, 
\begin{equation*} 	
\mathbb{E}_{\rm qs}\left[ X_r^n \right] = \pi_r n! - \theta \pi_r(1 - P_{rr}) ( n^2 - n -1 + n  H_n ) (n - 1)! + o(\theta);  
\end{equation*}
for $n_r, n_s > 0$, where $r \ne s \in \{1, \ldots, d\}$,   
\begin{equation*}  	
\mathbb{E}_{\rm qs}\left[ X_r^{n_r} X_s^{n_s} \right] = \theta \pi_r P_{rs} \frac{n_r!(n_s - 1)!}{n_s + 1}(n_r + 2n_s + 1)  + (r \leftrightarrow s) + o(\theta);  
\end{equation*}
and if three or more of the components of $(n_i, \ldots, n_d)$ are non-zero, 
\begin{equation*}
\mathbb{E}_{\rm qs}\left[ \prod_{i = 1}^d X_i^{n_i} \right] = o(\theta).
\end{equation*}
\end{theorem}
\begin{proof}
Consider first 
\begin{equation*}
\mathbb{E}_{\rm qs}\left[ X_r^n \right] = \int_0^\infty x_r^n \left( \sum_{j \ne r} \int_0^\infty g^{\rm surface} (x_r, x_j) dx_j + g^{\rm line} (x_r) \right) dx_r, 
\end{equation*}
with $g^{\rm surface}$ and $g^{\rm line}$ as given in Theorem~\ref{theorem2}.  The required integrals can be calculated using the following identities: 
\begin{equation*} \label{momentIdentities}
\begin{split}
&\int_0^\infty x^{a\theta + n -1}E_2(x) dx = \begin{cases}
	\displaystyle \frac{1}{a\theta} - (1 + \gamma) + \mathcal{O}(\theta) &n = 0; \\ \\
	\displaystyle \frac{(n - 1)!}{n + 1} + \mathcal{O}(\theta)  &  n \ge 1, 
	\end{cases} \\
&\int_0^\infty x^nE_1(x) dx = \frac{n!}{n + 1} \qquad n \ge 0, \\
&\int_0^\infty x^n e^{-x} dx = n!  \qquad n \ge 0, \\
&\int_0^\infty x^n e^{-x} \log x\, dx =\left(-\gamma + H_n \right) n!  \qquad n \ge 1.  \end{split}
\end{equation*}
The last identity in this list is a consequence of Gradshteyn and Ryzhik~\cite[][Eq.~(4.352.4)]{Gradshteyn:1965} and Abramowitz and Stegun~\cite[][Eqs.~(6.3.1) and (6.3.2)]{Abramowitz:1965sf}.  The $g^{\rm surface}$ integral 
contributes four parts: 
\begin{eqnarray}	\label{partA}
\lefteqn{ \int_0^\infty x_r^n \sum_{j \ne r} \int_0^\infty \theta \pi_j P_{jr} x_j e^{-x_j} x_r^{a_{rj}(x_j)\theta - 1} E_2(x_r) dx_j dx_r } \qquad\qquad \qquad\qquad \qquad\qquad\qquad\qquad \nonumber \\
	 & = & \theta \pi_r(1 - P_{rr}) \frac{(n - 1)!}{n + 1}; 
\end{eqnarray}
\begin{eqnarray}	\label{partB}
\lefteqn{ \int_0^\infty x_r^n \sum_{i \ne r} \int_0^\infty \theta \pi_r P_{ri} x_r e^{-x_r} x_i^{a_{ir}(x_r)\theta - 1} E_2(x_i) dx_i dx_r } \qquad\qquad \qquad\qquad\qquad\qquad \nonumber \\
	 & = & \pi_r n! - \theta \pi_r(1 - P_{rr}) (1 + \gamma) (n + 1)! ; 
\end{eqnarray}
\begin{eqnarray}	\label{partC}
\lefteqn{ \int_0^\infty x_r^n \sum_{j \ne r} \int_0^\infty \theta \pi_j P_{jr} 2e^{-x_j} E_1(x_r) dx_j dx_r } \qquad\qquad \qquad\qquad \qquad\qquad\qquad\qquad \nonumber \\
	 & = & 2\theta \pi_r(1 - P_{rr}) \frac{n!}{n + 1}; 
\end{eqnarray}
and 
\begin{eqnarray}	\label{partD}
\lefteqn{ \int_0^\infty x_r^n \sum_{i \ne r} \int_0^\infty \theta \pi_r P_{ri} 2e^{-x_r} E_1(x_i) dx_i dx_r } \qquad\qquad \qquad\qquad \qquad\qquad\qquad\qquad \nonumber \\
	 & = & 2\theta \pi_r(1 - P_{rr}) n! .
\end{eqnarray}
The line integral contributes a part 
\begin{eqnarray}	\label{partE}
\lefteqn{ -\theta \pi_r(1 - P_{rr}) \int_0^\infty x_r^n \left\{E_1(x_r) + [\gamma (1 - x_r) + \log x_r] e^{-x_r} \right\} dx_r} \qquad\qquad\qquad\qquad \nonumber \\
	& = & -  \theta \pi_r(1 - P_{rr})[H_{n + 1} - \gamma (n + 1) ] n! .
\end{eqnarray}
 
Adding Eqs.(\ref{partA}), (\ref{partB}), (\ref{partC}), (\ref{partD}) and (\ref{partE}) and simplifying gives, to $\mathcal{O}(\theta)$,  
\begin{equation*} 
\mathbb{E}_{\rm qs}\left[ X_r^n \right] = \pi_r n! - \theta \pi_r(1 - P_{rr}) ( n^2 - n -1 + n  H_n ) (n - 1)!  
\end{equation*}
as required.

Second, consider the case where $n_r, n_s > 0$ with $r \ne s$.  Then to $\mathcal{O}(\theta)$, 
\begin{eqnarray*} 
\mathbb{E}_{\rm qs}\left[ X_r^{n_r} X_s^{n_s} \right] & = & \int_0^\infty \int_0^\infty x_r^{n_r} x_s^{n_s} g^{\rm surface} (x_r, x_s) dx_r dx_s \nonumber \\
	& = & \int_0^\infty \int_0^\infty\theta \pi_s P_{sr} \left\{ x_s^{n_s + 1} e^{-x_s} x_r^{a_{rs}(x_s)\theta + n_r - 1} E_2(x_r) \right. \nonumber \\
	& & \qquad\qquad\qquad \left. \phantom{x_r^()}+ \, 2x_s^{n_s}e^{-x_s} x_r^{n_r} E_1(x_r) \right\} dx_r dx_s + (r \leftrightarrow s)  \nonumber \\
	& = & \theta \pi_s P_{sr} \left\{ \frac{(n_s + 1)!(n_r - 1)!}{n_r + 1} + \frac{2n_s! n_r!}{n_r + 1} \right\} + (r \leftrightarrow s)  \nonumber \\
	& = & \theta \pi_r P_{rs} \frac{n_r!(n_s - 1)!}{n_s + 1}(n_r + 2n_s + 1)  + (r \leftrightarrow s), 
\end{eqnarray*}
as required.

Clearly the presence of delta-functions in Eq.~(\ref{gToOrderTheta}) ensures that moments calculated to $\mathcal{O}(\theta)$ 
are identically zero if three or more of the components of $(n_1, \ldots, n_d)$ are non-zero.  
\end{proof}

%
\subsection{Moments in $\mathbf{U} = \mathbf{X}/X$  to order $\theta$ and sampling distributions}

We are also interested in moments of the relative proportions of each allele type, as this will enable calculation of sampling distributions.  
%
%
\begin{theorem}\label{theorem4}
Define the total population and relative proportion of each allele type respectively as 
\begin{equation*}
X = \sum_{i = 1}^d X_i, \qquad U_i = \frac{X _i}{X}, \quad i = 1, \ldots d, 
\end{equation*}
where only $d - 1$ of the $U_i$ are independent because of the constraint $\sum_{i = 1}^d U_i = 1$.  
Then we have the following moments for the asymptotic relative proportions: For integer $n \ge 1$ and $r \in \{1, \ldots, d\}$,   
\begin{equation*}
\mathbb{E}_{\rm qs}\left[ U_r^n \right] = \pi_r \left(1 - \theta(1 - P_{rr}) H_{n - 1} \right) + o(\theta);  
\end{equation*}
for $n_r, n_s > 0$, where $r \ne s \in \{1, \ldots, d\}$, 
\begin{equation*}
\mathbb{E}_{\rm qs}\left[ U_r^{n_r} U_s^{n_s} \right] = 
				\theta\pi_s P_{sr} \frac{(n_r - 1)! n_s!}{(n_r + n_s)!}  + (r \leftrightarrow s) + o(\theta); 
\end{equation*}
and if three or more of the components of $(n_i, \ldots, n_d)$ are non-zero, 
\begin{equation*}
\mathbb{E}_{\rm qs}\left[ \prod_{i = 1}^d U_i^{n_i} \right] = o(\theta).
\end{equation*}
\end{theorem}
\begin{proof}
The density of the quasi-stationary distribution corresponding to the random variables $(X, \mathbf{U})$ is, from Eq.~(\ref{gToOrderTheta}) 
and the fact that $\delta(xu_\ell) = x^{-1}\delta(u_\ell)$,
\begin{eqnarray*}
\lefteqn{g_{X, \mathbf{U}}(x, \mathbf{u}) = x^{d - 1} g_{\mathbf{X}}(x \mathbf{u}) \delta\left(1 - \sum_{i = 1}^d u_i \right)} \nonumber \\
	& = &  \left(\sum_{1 \le i < j \le d} x g^{\rm surface} (xu_i, xu_j) \prod_{\ell \ne i, j} \delta(u_\ell) + 
											 \sum_{i = 1}^d g^{\rm line} (xu_i) \prod_{\ell \ne i} \delta(u_\ell)\right) \nonumber \\ 
	& & \qquad\qquad\qquad\qquad\qquad\qquad\qquad\qquad\qquad \times \delta\left(1 - \sum_{i = 1}^d u_i \right).
\end{eqnarray*}
Once again, the presence of delta-functions ensures that moments calculated to $\mathcal{O}(\theta)$ are identically zero if 
three or more of the components of $\mathbf{n} = (n_1, \ldots, n_d)$ are non-zero.  Thus only two cases need be considered.  

For $n_r, n_s > 0$ where $r \ne s \in \{1, \ldots, d\}$, 
\begin{eqnarray}	\label{UnrUnsMoments}
\mathbb{E}_{\rm qs}\left[ U_r^{n_r} U_s^{n_s} \right] 
	& = & \int_0^1 \cdots \int_0^1 u_r^{n_r} u_s^{n_s} \int_0^\infty g_{X, \mathbf{U}}(x, \mathbf{u}) dx\,  du_1 \cdots du_d \nonumber \\
	& = &  \int_0^1  \int_0^1 u_r^{n_r} u_s^{n_s} \int_0^\infty x g^{\rm surface} (xu_r, xu_s) \delta(1 - u_r - u_s) dx\,  du_r du_s \nonumber\\
	& = & \theta\pi_s P_{sr} \int_0^1  u^{n_r} (1 - u)^{n_s} (\mathcal{I}_1(u) + \mathcal{I}_2(u))  du + ( r \leftrightarrow s), 
\end{eqnarray} 
where 
\begin{eqnarray*}
\mathcal{I}_1(u) & = & \int_0^\infty 2xe^{-(1 - u)x} E_1(ux) dx, \nonumber \\
	& = & - \frac{2}{1 - u} - \frac{2 \log u}{(1 - u)^2}, 
\end{eqnarray*}
\begin{eqnarray*}
\mathcal{I}_2(u) & = & \int_0^\infty \frac{1 - u}{u} x e^{-(1 - u)x} E_2(ux) dx \nonumber \\
	& = & \frac{1}{u} + \frac{2}{1 - u} + \frac{2 \log u}{(1 - u)^2}.   
\end{eqnarray*}
Note that in Eq.~(\ref{UnrUnsMoments}) 
the factors $u_r^{n_r} u_s^{n_s}$ in combination with the $\prod_\ell \delta(u_\ell)$ ensure that the surface terms but not the line densities survive the integration,
and that to first order 
in $\theta$ the $a_{ij}$ term in the exponent of Eq.~(\ref{gSurface}) can be set to zero provided $n_r, n_s > 0$.  The last two integrals have 
been evaluated using Wolfram$|$Alpha~\citep{WolframAlpha}
with the code
\begin{quote}
\tt{Integrate(2*x*Exp(-(1 - u)*x) * ExpIntegral[1, u*x]) from x=0 to x=infinity} 
\end{quote} 
and
\begin{quote}
\tt{Integrate((1 - u)/u * x * Exp(-(1 - u)*x) * ExpIntegral[2, u*x]) from x=0 to x=infinity} 
\end{quote} 
respectively.  Then 
\begin{eqnarray*}
\mathbb{E}_{\rm qs}\left[ U_r^{n_r} U_s^{n_s} \right] & = & \theta\pi_s P_{sr} B(n_r, n_s + 1) + (r \leftrightarrow s) \nonumber \\
	& = &\theta\pi_s P_{sr} \frac{(n_r - 1)! n_s!}{(n_r + n_s)!}  + (r \leftrightarrow s). 
\end{eqnarray*}

For $n > 1$, the above result leads to the iterative rule 
\begin{eqnarray*}
\mathbb{E}_{\rm qs}\left[ U_r^n \right] & = & \mathbb{E}_{\rm qs}\left[ U_r^{n - 1} \left(1 - \sum_{\ell \ne r} U_\ell \right) \right] \nonumber \\
	& = & \mathbb{E}_{\rm qs}\left[ U_r^{n - 1} \right] - \theta \pi_r(1 - P_{rr}) \frac{1}{n - 1}. 
\end{eqnarray*}
For $n = 1$, the asymptotic probability of observing a single individual sampled from a surviving population to be of type-$r$ is 
$\mathbb{E}_{\rm qs}\left[ U_r \right] = \pi_r$.  Hence 
\begin{equation*}
\mathbb{E}_{\rm qs}\left[ U_r^n \right] = \pi_r \left(1 - \theta(1 - P_{rr}) H_{n - 1} \right).  
\end{equation*}
\end{proof}
%
%
\begin{corollary}
In a random sample of $n$ individuals from 
the quasi-stationary limit of a subcritical multi-type branching diffusion, the probability that the $d$ types are distributed within the sample
as $\mathbf{n} = (n_1, \ldots, n_d)$, where $\sum_{i = 1}^d n_i = n$, is
\begin{equation}	\label{samplingDistrib}
\begin{split}
p(n \mathbf{e}_r) & = \pi_r \left(1 - \theta(1 - P_{rr}) H_{n - 1} \right) + \mathcal{O}(\theta^2) \\
p(n_r \mathbf{e}_r + n_s \mathbf{e}_s) &= \theta\left(\pi_r P_{rs} \frac{1}{n_s} + \pi_s P_{sr} \frac{1}{n_r} \right)  + \mathcal{O}(\theta^2) \\
p(\mathbf{n}) &=  \mathcal{O}(\theta^2) \quad\text{if $\mathbf{n}$ has $>2$ non-zero entries,}
\end{split}
\end{equation}
as $\theta \to 0$, where $r \ne s \in \{1, \ldots d\}$ and $n_r + n_s = n$. 
\end{corollary}
\begin{proof}
The required sampling distribution is an immediate consequence of Theorem~\ref{theorem4} and the formula 
\begin{equation*}
p(\mathbf{n}) = {n \choose \mathbf{n}} \mathbb{E}_{\rm qs}\left[ U_1^{n_1} \cdots U_d^{n_d} \right]. 
\end{equation*}
\end{proof}

This distribution is identical to the multi-allele stationary sampling distribution for a neutral Wright-Fisher 
population, first determined by \citet{BurdenTang17}, and subsequently verified using alternate methods by 
\citet{SchrempfHobolth17} and \citet{Burden_2018}.  

Finally, recall from Subsection~\ref{sec:SmallRates} that the parameter $\alpha$ can be reinstated by making the substitution 
$\theta \to \tfrac{1}{2}|\alpha|^{-1}\theta$ where $\alpha$ is related to the discrete BGW branching process by Eq.~(\ref{alphaDef}), 
and it is this combination which is assumed to be small.  
%
%

\section{Comparison with numerical simulation for $d=2$ types}
\label{sec:ComparisonWithNumerical}

Suppose the population $Y(\tau)$ of the discrete BGW process described in Section~\ref{sec:MultitypeBranching} is divided into $d = 2$ types of size $Y_1(\tau)$ 
and $Y_2(\tau) = Y(\tau) - Y_1(\tau)$ respectively, with per-generation mutation rates between the two types $r_{12}$ and $r_{21}$.   
Define a transition probability 
\begin{eqnarray}	\label{discreteMultitypeProcess} 
P(m, i; n, j) & := & \Prob(Y(\tau + 1) = n, Y_1(\tau + 1) = j \mid Y(\tau) = m, Y_1(\tau) =i) \nonumber \\
	& = & p(m, n) {n \choose j} \chi(i, m)^j (1 - \chi(i, m))^{n - j}, 
\end{eqnarray}
for $m, n = 0, 1, 2, \ldots$; $i = 0, \ldots, m$ and $j = 0, \ldots, n$, where 
\begin{equation*}
\chi(i, m) = \frac{i}{m}(1 - r_{12}) + \left(1 - \frac{i}{m}\right) r_{21}, 
\end{equation*}
and $p(m, n) = \Prob(Y(\tau + 1) = n \mid Y(\tau) = m)$.  
Note that $Y(\tau) = 0$ is an absorbing state corresponding to extinction of the entire population.  
Our aim is to compare a numerical determination of the quasi-stationary distribution of this transition matrix with the theoretical small-rates continuum 
diffusion limit density derived in Section~\ref{sec:QSSubcriticalMultiType}.  

The quasi-stationary distribution, if it exists, will be of the form 
\begin{equation*}
G(m, i) = \begin{cases}
	0 & \text{if } m = 0, \\
	\tilde{G}(m, i) & \text{if } m > 0 \text{ and } i = 0, \ldots, m, 
	\end{cases}
\end{equation*}
where $\tilde{G}$ is a left eigenvector of a matrix $\tilde{P}$, equal to the transition matrix $P(m, i; n, j)$ with the first row ($n = 0$) and first column ($m = 0$) 
removed.   To see this, observe that updating $G$ by one time step results in the state $G P = (\Pi, (1 - \Pi) \tilde{G})$, where $\Pi$ is the limiting probability 
of extinction in one time step as $\tau \to \infty$ given survival of the population to time $\tau$.   Thus the quasi-stationary distribution is obtained numerically 
by computing the principal left-eigenvector of $\tilde{P}$ and renormalising the sum of the elements to 1. 
 
For $d=2$, this distribution is to be matched with the single 2-dimensional surface density $g^{\rm surface}(x_1, x_2)$ in Eq.~(\ref{gSurface}).  
The scaled populations in the diffusion limit corresponding to $Y_i(\tau)$ are found from Eqs.~(\ref{diffLimit}) and (\ref{alphaDef}) to be 
\begin{equation}	\label{YtoXscaling}
\mathbf{X}(t) \approx \frac{\log\lambda}{\alpha \sigma^2} \mathbf{Y}(\tau). 
\end{equation}
Then setting 
\begin{equation*}
(x_1, x_2) = \frac{\log\lambda}{\alpha \sigma^2} (i, m - i),  \quad dx_1 = dx_2 = \frac{\log\lambda}{\alpha \sigma^2}, 
\end{equation*}
and applying a coordinate transformation 
\begin{equation}	\label{XToXUTransf} 
\quad x = x_1 + x_2, \quad u = x_1/x, \quad g^{\rm surface}(x_1, x_2) = x^{-1} g^{\rm surface}_{X, U}(x, u), 
\end{equation}
we have 
\begin{eqnarray*}
\tilde{G}(m, i) & = & \lim_{\tau \to \infty}\Prob(Y_1(\tau) = i, Y_2(\tau) = m - i \mid Y(\tau) > 0) \nonumber \\
	& \approx & g^{\rm surface}(x_1, x_2) dx_1 dx_2 \nonumber \\
	& = & \frac{\log\lambda}{m\alpha \sigma^2} g^{\rm surface}_{X, U}(x, u), 
\end{eqnarray*}
or 
\begin{equation}	\label{discreteToSurfaceDensity}
g^{\rm surface}_{X, U}(x, u) \approx \frac{m\alpha \sigma^2}{\log\lambda} \tilde{G}(m, i). 
\end{equation}
The marginal probability in the total population size is related to the diffusion limit via $g_{X}(x)dx \approx \sum_{i = 0}^m\tilde{G}(m, i)$, or 
\begin{equation}	\label{discreteToMarginalX}
g_{X}(x) \approx \frac{\alpha \sigma^2}{\log\lambda} \sum_{i = 0}^m\tilde{G}(m, i). 
\end{equation}

\begin{figure}[t!]
\begin{center}
\centerline{\includegraphics[width=0.76\textwidth]{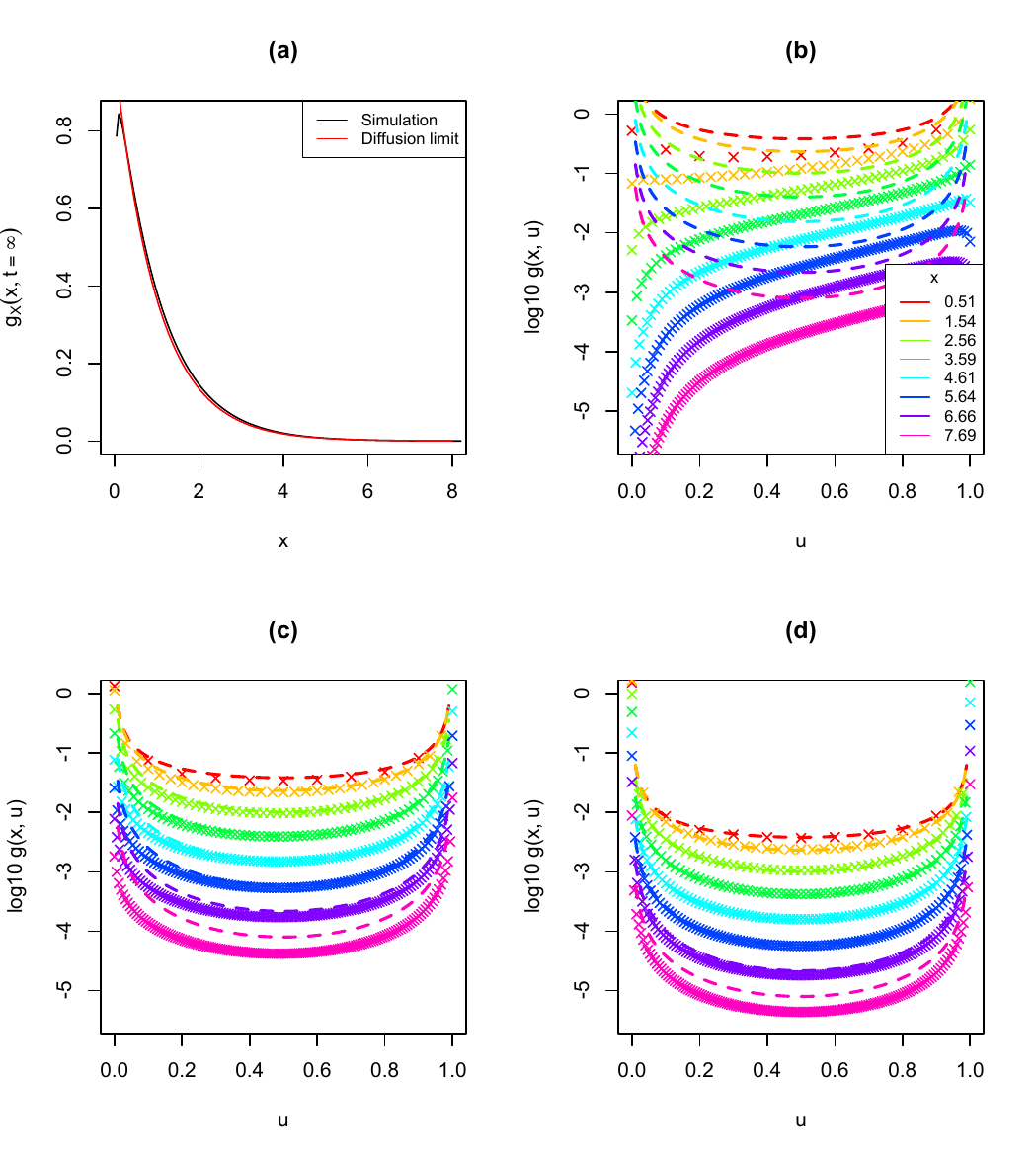}}
\caption{Comparison of the numerically determined quasi-stationary distribution of the subcritical branching process, Eq.~(\ref{discreteMultitypeProcess}), 
with the theoretical small-rates continuum diffusion limit density (see Section~\ref{sec:QSSubcriticalMultiType}):  (a) Comparison of the computed marginal 
distribution in the total population size scaled to a continuum density via Eq.~(\ref{discreteToMarginalX}) (black curve) with the diffusion limit exponential 
density Eq.~(\ref{gExponential}) (red curve); (b) to (d): comparison of the computed quasi-stationary distribution scaled to a continuum surface density via 
Eq.~(\ref{discreteToSurfaceDensity}) (crosses) with the surface density Eq.~(\ref{gSurface}) for $\theta = 1$, 0.1 and 0.01 respectively, and 
$(\pi_1, \pi_2) = (0.75, 0.25)$ (dashed lines).  The remaining parameters are set as described in the text.} 
\label{fig:SurfaceDensitySimulation}
\end{center}
\end{figure}

Figure~\ref{fig:SurfaceDensitySimulation} shows plots of the computed discrete quasi-stationary distribution transformed to a surface density 
via Eqs.~(\ref{discreteToMarginalX}) and (\ref{discreteToSurfaceDensity}).  Superimposed are plots of the theoretical densities Eq.~(\ref{gExponential}) and Eq.~(\ref{gSurface}) 
transformed to $(x, u)$ coordinates via Eq.~(\ref{XToXUTransf}).  For simplicity we choose the distribution of the number of offspring per parent to be Poisson, 
\begin{equation*}
P(m, n) = \frac{e^{-\lambda m}(\lambda m)^n}{n!},  
\end{equation*}
where $\log\lambda = -0.025$, and thus mean number of offspring per parent $\lambda = \sigma^2 \approx 0.9753$.  
To compute the principal eigenvector of $P(m, i; n, j)$ a cutoff $m_{\rm max} = 160$ is implemented 
on total population size.  With $\alpha = -\tfrac{1}{2}$, the corresponding cutoff on the diffusion limit population size is 
$x_{\rm max} = m_{\rm max} \log \lambda/(\alpha\sigma^2) \approx 8.2$, so that truncation of the exponential density Eq.~(\ref{gExponential}) removes a fraction 
no more than $e^{-8.2} \approx 0.00027$ of the total probability.  Figure~\ref{fig:SurfaceDensitySimulation}(a) compares the computed marginal quasi-stationary 
distribution scaled to a continuum density via Eq.~(\ref{discreteToMarginalX}) with the diffusion limit exponential density.  The close agreement confirms the 
suitability of the chosen parameters $\lambda$ and $m_{\rm max}$.   It remains to choose $r_{12}$ and $r_{21}$ in a way that will enable the range of validity of 
the small expansion parameter $\theta$ to be determined.  

Because of the scale invariance of the quasi-stationary distribution, Eq.(\ref{gScaling}), the only free parameters in the diffusion limit are the 
rate matrix elements $\gamma_{ij}$ relevant to $\alpha = -\tfrac{1}{2}$, and these can be specified in terms of the parameters $\theta$ and $P_{ij}$.  
From Eqs.~(\ref{alphaDef}) and (\ref{YtoXscaling}) the diffusion limit rate matrix is 
\begin{equation*}
\tfrac{1}{2}\theta P_{ij} = \gamma_{ij} = \frac{\alpha}{\log\lambda} r_{ij}, \qquad i \ne j.  
\end{equation*}
For $d=2$ the PIM form of the rate matrix,  $P_{ij} = \pi_j$ is appropriate, so 
\begin{equation*}
r_{ij} = \tfrac{1}{2} \theta \pi_j \frac{\log\lambda}{\alpha}.  
\end{equation*}
Figures~\ref{fig:SurfaceDensitySimulation}(b) to (d) compare the computed quasi-stationary distribution scaled to a continuum surface density via 
Eq.~(\ref{discreteToSurfaceDensity}) with the surface density Eq.~(\ref{gSurface}) for $\theta = 1$, 0.1 and 0.01 respectively, 
and $(\pi_1, \pi_2) = (0.75, 0.25)$.  One see that the small-rates approximation to first order in $\theta$ performs well provided $\theta \le 0.1$, and 
poorly for $\theta$ of order unity.  
Disagreement between simulation and theory at values of $x$ approaching the cutoff at $x_{\rm max} = 8.2$ in 
Figs.~\ref{fig:SurfaceDensitySimulation}(c) and (d) is mainly due to comparing a discrete distribution with an imposed hard cutoff 
on the total population size with the infinite tail of the diffusion distribution.  The difference is amplified by the logarithmic scale of the plot.

\section{Conclusions}
\label{sec:Conclusions}

Certain asymptotic properties of discrete multi-type branching processes have been well known for some time~\citep{Harris64}.  Here we have approached 
the topic directly from the continuum viewpoint of the diffusion limit.  There are two advantages to this approach.  Firstly, the approach is accessible to population 
geneticists, who are well aware of the influence of Kimura's use of forward Kolmogorov equations to study the fixation of allelic mutations in populations.  
Use of the diffusion limit in population genetics is dominated by Wright-Fisher, Moran or similar models constrained so that total population size is set externally.  There 
have been relatively few treatments in the applied population genetics literature acknowledging a population whose size is determined stochastically.  
Secondly, from a mathematical point of view, some results can be more readily obtained from the diffusion process than from the discrete process.  

Our treatment has concentrated on neutral mutations.  This enables us to exploit the mathematical simplification that the total population size is effectively a 
Feller diffusion for a single allele type.  For subcritical and critical process the population goes extinct almost surely and the interesting limit is the Yaglom limit conditioned 
on non-extinction~\citep{yaglom1947certain}.  

Our calculation of the stationary properties of the supercritical and critical multi-type diffusions in Section~\ref{sec:AsymptoticSuperAndCrit} aim to provide 
easily accessible derivations of known results in the formal continuous-state branching process literature~\cite{champagnat2008limit,Kyprianou18}
The calculation of the stationary distribution in the supercritical case generalises an earlier result for 2 types to the general case of $d$ types~\citep{burden2018mutation}.  
The resulting distribution with the exponential growth factored out is the analogue of the known result for a discrete branching process, namely a one-dimensional line 
density directed along a ray aligned with the stationary eigenvector of the rate matrix.  The line density in the continuum limit is equal to 
the solution by \citet{feller1951diffusion}.  A similar result follows for the quasi-stationary critical case, except that the relevant density is that of the population with linear time 
factored out, and the asymptotic line density agrees with Yaglom's exponential distribution.  

The main results of this paper in Sections~\ref{sec:QSSubcriticalMultiType} to \ref{sec:ComparisonWithNumerical} pertain to the subcritical branching diffusion, 
for which the quasi-stationary distribution does not collapse on to a line density.  
Although an exact quasi-stationary distribution remains intractable, 
a solution is found to first order in the overall mutation rate $\theta$ via a multi-dimensional Laplace transform leading to a first-order partial differential 
equation, which we solve using the method of characteristics.  The solution agrees well with numerically determined quasi-stationary distributions of discrete 
multi-type branching processes provided $\theta < 0.1$.  As an order of magnitude estimate, $\theta$ can be thought of as the product of a per base 
mutation rate per nucleotide site per generation and an effective population size, and this product is less than 0.1 in most biological 
contexts~\citep[see][Fig.~3b]{Lynch16}. 

Of particular interest is our calculation from the marginal distribution of the relative proportion of allele types of 
the sampling distribution over types for a sample of given finite size (see Eq.~(\ref{samplingDistrib})).  This sampling distribution is identical to 
the $\mathcal{O}(\theta)$ multi-allele stationary sampling distribution for a neutral Wright-Fisher diffusion with fixed population size~\citep{Burden_2018}.  

%
%
\section*{Declaration of interest statement}

No potential competing interest was reported by the authors.

\bibliographystyle{plainnat}

\end{document}